\documentclass{amsart}
\usepackage{amssymb}
\usepackage{mathrsfs}
\usepackage{amsrefs}
\usepackage{array}
\usepackage{ulem}

\usepackage[all]{xy}
\usepackage{tikz}

\usepackage{multicol}

\usepackage{hyperref}  
\hypersetup{%
  bookmarksnumbered=true,%
  bookmarks=true,%
  colorlinks=true,%
  linkcolor=blue,%
  citecolor=blue,%
  filecolor=blue,%
  menucolor=blue,%
  pagecolor=blue,%
  urlcolor=blue,%
  pdfnewwindow=true,%
  pdfstartview=FitBH}

\newcommand{\al}{\alpha}
\newcommand{\be}{\beta}
\newcommand{\ga}{\gamma}

\newcommand{\SI}{\Sigma}

\newcommand{\F}{\mathbb{F}_2}

\newcommand{\R}{\mathbb{R}}
\newcommand{\Z}{\mathbb{Z}}

\newcommand{\xrtar}{\xrightarrow}
\newcommand{\xrtarr}[1]{\stackrel{#1}{\longrightarrow}}
\newcommand{\rtarr}{\longrightarrow}
\newcommand{\into}{\hookrightarrow}
\newcommand{\onto}{\twoheadrightarrow}
\newcommand{\iso}{\cong}
\newcommand{\sma}{\wedge}
\DeclareMathOperator{\fib}{fib}
\newcommand{\ul}{\underline}
\newcommand{\ulF}{\ul{\F}}
\newcommand{\ulf}{\ul{f}}
\newcommand{\ulg}{\ul{g}}
\newcommand{\mf}[1]{\ul{#1}}
\newcommand{\mpi}{\ul{\pi}}
\DeclareMathOperator{\Hom}{Hom}
\DeclareMathOperator{\coker}{coker}

\numberwithin{equation}{section}
\numberwithin{figure}{section}

\def\makeautorefname#1#2{\expandafter\def\csname#1autorefname\endcsname{#2}}
\makeautorefname{theorem}{Theorem}%
\makeautorefname{thm}{Theorem}%
\makeautorefname{addm}{Addendum}%
\makeautorefname{mainthm}{Main theorem}%
\makeautorefname{corollary}{Corollary}%
\makeautorefname{const}{Construction}%
\makeautorefname{cor}{Corollary}%
\makeautorefname{lemma}{Lemma}%
\makeautorefname{notn}{Notation}%
\makeautorefname{lem}{Lemma}%
\makeautorefname{sublemma}{Sublemma}%
\makeautorefname{sublem}{Sublemma}%
\makeautorefname{subl}{Sublemma}%
\makeautorefname{prop}{Proposition}%
\makeautorefname{warn}{Warning}%
\makeautorefname{eg}{Example}%
\makeautorefname{figure}{Figure}%

\newtheorem{thm}{Theorem}[section]
\newtheorem{cor}{Corollary}[section]
\newtheorem{prop}{Proposition}[section]
\newtheorem{lem}{Lemma}[section]

\newtheorem*{mainres}{Main Result}

\theoremstyle{definition}
\newtheorem{defn}{Definition}[section]

\newtheorem{eg}{Example}[section]
\newtheorem{notn}{Notation}[section]
\newtheorem{rem}{Remark}[section]

\newenvironment{pf}{\begin{proof}}{\end{proof}}

\makeatletter
\let\c@cor=\c@thm
\let\c@prop=\c@thm
\let\c@lem=\c@thm
\let\c@defn=\c@thm
\let\c@eg=\c@thm
\let\c@notn=\c@thm
\let\c@rem=\c@thm
\let\c@warn=\c@thm
\let\c@sch=\c@thm
\let\c@equation=\c@thm
\let\c@figure=\c@thm
\makeatother

\newcolumntype{L}{>{$}l<{$}}

\newenvironment{psmallmatrix}
  {\left(\begin{smallmatrix}}
  {\end{smallmatrix}\right)}
\newcommand{\smallmtx}[1]{\begin{psmallmatrix}#1\end{psmallmatrix}}

\newcounter{themyfigure}

\newcommand{\Kl}{\mathcal{K}}

\newcommand{\MackC}[2]
{\xymatrix{
 #1 \ar@/_2ex/[d]  \\
 #2  \ar@/_2ex/[u] 
} }
\newcommand{\MackCAr}[4]
{\xymatrix{
 #1 \ar@/_2ex/[d]_{#3}  \\
 #2  \ar@/_2ex/[u]_{#4} 
} }
\newcommand{\MackK}[5]
{\xymatrix{
 & #1 \ar@/_1ex/[dl] \ar@/_1ex/[d] \ar@/_1ex/[dr]& \\
 #2 \ar@/_1ex/[dr] \ar@/_1ex/[ur] & #3 \ar@/_1ex/[d] \ar@/_1ex/[u]& #4 \ar@/_1ex/[ul]  \ar@/_1.5ex/[dl] \\
  & #5 \ar@/_1ex/[ul] \ar@/_1ex/[u] \ar@/_1.5ex/[ur]&
} }

\title{The Klein four slices of $\SI^n H\protect\underline{\mathbb{F}}_2$}
\author{B. Guillou}
\address{Department of Mathematics, The University of Kentucky, Lexington, KY 40506--0027}
\email{bertguillou@uky.edu}
\author{C. Yarnall}
\address{Department of Mathematics, California State University -
  Dominguez Hills, Carson, CA, 90747}
\email{cyarnall@csudh.edu}

\thanks{B. Guillou was  supported by NSF grant DMS-1710379.}

\subjclass{Primary 55N91,  55P91; 55Q91 \\Secondary 55T99}

\begin{document}

\begin{abstract}
We describe the slices of positive integral suspensions of the equivariant Eilenberg-Mac~Lane spectrum $H\ulF$ for the constant Mackey functor 
over the Klein four-group $C_2\times C_2$.
\end{abstract}

\maketitle
\tableofcontents

\section{Introduction}
The slice filtration is a filtration of equivariant spectra developed by Hill, Hopkins, and Ravenel, as a
generalization of Dugger's filtration \cite{D}, in
their solution to the Kervaire invariant-one problem \cite{Kervaire}. It is
an equivariant analogue of the Postnikov tower and was
modeled on the motivic filtration of Voevodsky \cite{V}. 

Since its inception, there have been a few reformulations and new
understandings of the structure of the slice filtration. Some
properties and useful results in this setting are summarized in
\autoref{sec:slicefiltration}. In this paper, we use the {\bf regular} slice filtration (cf. \cite{U}, \cite{HY}) on
equivariant spectra and note that this filtration differs from the
original filtration from \cite{Kervaire} by a shift by one. 

Let $G$ be a finite group and let $\mathrm{Sp}^G$ be the category of genuine $G$-spectra. 

\begin{defn}
Let $\tau_{\geq n}^G\subseteq \mathrm{Sp}^G$ be the localizing
subcategory generated by $G$-spectra of the form $\Sigma^\infty_G G/H_+
\sma S^{k \rho_H}$, where $H \subset G$, $\rho_H$ is the regular representation of
$H$ and $k\cdot |H| \geq n$.
We  write $X\geq n$ to mean that $X\in \tau_{\geq n}^G$.
\end{defn}

We use $P^{n-1}(-)$ to denote the localization functor associated to
$\tau_{\geq n}^G$. There are natural transformations $P^n(-) \rtarr
P^{n-1}(-)$ that give the slice tower of $X$
\[ \cdots \rtarr P^{n+1}X \rtarr P^n X \rtarr P^{n-1} X \rtarr \dots,\]
and the fiber at each level
\[ P^n_n X \rtarr P^n X \rtarr P^{n-1}X\]
is known as the {\bf $n$-slice} of $X$.

While in the nonequivariant setting the relationship between the
Postnikov tower of a spectrum and the spectrum's homotopy groups is clear, there
is a much more complicated story for homotopy groups and the slice tower when working
equivariantly. Furthermore, such homotopy groups enjoy a richer structure. For a $G$-spectrum $X$, the homotopy groups
$\pi_n(X^H)$, as $H$ varies over the subgroups of $G$, define a
$G$-Mackey functor. An underline will denote a Mackey functor, and we
will display such functors $\mf{M}$ according to
their Lewis diagrams. The general form of such diagrams for $G = C_2$ and
$G = C_2 \times C_2$ are displayed below\footnote{We write $\mf{M}(L)$ for what would be typically written as $\mf{M}((C_2\times C_2)/L)$}.
\begin{center}
\[
\raisebox{-3ex}{\MackC{\mf{M}(C_2)}{\mf{M}(e)}} \hspace{2cm}
\MackK{\mf{M}(C_2\times C_2)}{\mf{M}(L)}{\mf{M}(D)}{\mf{M}(R)}{\mf{M}(e)}\]
\end{center}

Here $L$, $D$, and $R$ are the left, diagonal, and right cyclic subgroups of $C_2\times C_2$ of  order two. 
We have not drawn in the Weyl group actions
on the intermediate groups or the $G$-action on $\mf{M}(e)$.
The maps pointing down are called restriction, and the maps pointing up are called transfers.

Associated to every $G$-Mackey functor $\mf{M}$, there is an
Eilenberg-MacLane spectrum $H\mf{M}$. While $H\mf{M}$ is always a
$0$-slice, and thus has a trivial slice tower, suspensions of
Eilenberg-MacLane spectra produce interesting slices and corresponding towers. For
instance, when $G = C_{p^n}$ the slices of $\Sigma^n H\mf{\Z}$ and
$\Sigma^{n\lambda}H\mf{\Z}$, where $\mf{\Z}$ is the constant Mackey
functor at $\Z$ and $\lambda$ is an irreducible $C_{p^n}$-representation, were presented in \cite{Y} and \cite{HHR}, respectively.

We will primarily work with the constant functor $\ulF$ for the Klein four-group $C_2\times C_2$ (which we will often denote by $\Kl$).
This Mackey functor takes on the value $\F$ at each subgroup. The restriction maps are all the identity, and the transfers are all zero.
In this paper we present the slices of  $\SI^n H\ulF$ for the Klein
four-group $\Kl=C_2\times C_2$ and $n\geq 0$. A summary of our main
results is as follows:

\begin{mainres} For $\Kl = C_2 \times C_2$ and $n \geq 0$, all
nontrivial slices of $X=\SI^n H_\Kl \ulF$ are given by:
\[ P^i_i(X) = \SI^V H_\Kl \mf{M}\]
where $|V| = i$ and $i \equiv 0 \mbox{ (mod }4)$ for $ n\leq i \leq 4n-12$, $i \equiv 2 \mbox{
  (mod } 4)$ for $n \leq i \leq 2n-4$, or $i = n$. The precise representations $V$ and Mackey
functors $\mf{M}$ are completely described in \autoref{nSliceProp},
\autoref{4kSlicesThm}, and \autoref{4k2SlicesProp}.
\end{mainres}

The paper is organized as follows. We begin with some background material in \autoref{sec:background}. In \autoref{sec:C2}, we review results from \cite{Kervaire} for the case of $H_{C_2}\ulF$. We present the relevant $\Kl$-Mackey functors  in \autoref{sec:KMackey}. Our main results, which describe all of the slices of $\SI^n H_\Kl \ulF$ are given in \autoref{sec:slices}. In \autoref{sec:towers}, we present the first few slice towers (up to $n=8$). The homotopy Mackey functors of the slices are computed in \autoref{sec:Homotopy}. Finally, in \autoref{sec:sliceSS}, we display a few examples of the slice spectral sequence for $\SI^n H_\Kl \ulF$. For convenience, we also list the important $\Kl$-Mackey functors in the \nameref{sec:app}.

We are grateful to John Greenlees, Mike Hill, Doug Ravenel, Nicolas Ricka, and Dylan Wilson for some helpful conversations.
We would like to thank the Isaac Newton Institute for Mathematical Sciences, Cambridge, for support and hospitality during the programme ``Equivariant and motivic homotopy theory'', where work on this paper was completed. Figures \ref{78Fig}, \ref{910Fig}, \ref{1112Fig}, and \ref{20Fig} were created using Hood Chatham's {\tt spectralsequences} package.

\section{Background}\label{sec:background}

\subsection{$(C_2 \times C_2)$-representations}
Recall that the real representation ring for the group $\Kl=C_2\times C_2$ is 
\[RO(\Kl) \iso \Z\{1, \alpha_{1,0}, \alpha_{1,1}, \alpha_{0,1}\}, \]
 where $1$ is the trivial one-dimensional representation and the other representations are defined by
\[\begin{split} \Z/2 \times \Z/2 &\xrtar{\alpha_{i,j}} \Z/2\stackrel{\sigma}{\into} \mathrm{Gl}_1(\R) \\
(k,n) & \mapsto ik+jn.
\end{split} \]
Thus $\alpha_{1,0}$ is the projection onto the left factor. 
To avoid cluttering notation, we prefer to write $\alpha = \alpha_{1,0}$, $\beta = \alpha_{0,1}$, $\gamma = \alpha_{1,1}$. We denote by $\rho$ or $\rho_\Kl$ the regular representation, and we have
\[ \rho = 1 + \alpha + \beta + \gamma\]
in $RO(\Kl)$.
The left, diagonal, and right cyclic subgroups of $\Kl$ will be denoted by $L$, $D$, and $R$, respectively. We have
\[ L = \ker\beta, \qquad D=\ker\gamma, \qquad R=\ker\alpha.\]

It will often be important to consider restriction to the cyclic subgroups. Given that $RO(C_2) \iso \Z\{1,\sigma\}$, the restrictions of representations are given by 
\[ RO(\Kl) \xrtarr{\iota^*} RO(C_2),\]
\[ \iota^*_L = \smallmtx{1 & 0 & 0 & 1 \\ 0 & 1 & 1 & 0 }, \qquad
\iota^*_D = \smallmtx{1 & 0 & 1 & 0 \\ 0 & 1 & 0 & 1 }, \qquad
\iota^*_R = \smallmtx{1 & 1 & 0 & 0 \\ 0 & 0 & 1 & 1 }
\]
In particular, we have $\iota^*(\rho_\Kl) = 2\rho_{C_2}$ in $RO(C_2)$.

Since the subgroups $H=L,D,R\unlhd \Kl$ are all normal, we get an induced action of $C_2\iso \Kl/H$ on the $H$-fixed points of any $\Kl$-representation.
These fixed point homomorphisms are given by
\[ RO(\Kl) \xrtar{(-)^H} RO(C_2),\]
\[ (-)^L = \smallmtx{1 & 0 & 0 & 0 \\ 0 & 0 & 0 & 1 }, \qquad
(-)^D = \smallmtx{1 & 0 & 0 & 0 \\ 0 & 0 & 1 & 0  }, \qquad
(-)^R = \smallmtx{1 & 0 & 0 & 0 \\ 0 & 1 & 0 & 0  }
\]
In particular, for any of these index two subgroups $H \leq \Kl$, we
have $(\rho_\Kl)^H = \rho_{C_2}$.

\subsection{Mackey functors}

For $G$-spectra $W$ and $X$, the collection of abelian groups $[ G/H_+\sma W ,X]^G$, as $H$ varies, defines a $G$-Mackey functor. In the case $W=S^V$ for a (virtual) $G$-representation $V$, this is the Mackey functor $\mpi_V(X)$. We give examples of $G$-Mackey functors for $G=C_2$ in \autoref{sec:C2Mackey} and for $G=C_2\times C_2$ in \autoref{sec:KMackey}.

\begin{notn}
We will typically denote Mackey functor restriction maps by 
\[r\!\downarrow^G_{H}:\mf{M}(G)\rtarr \mf{M}(H)\]
 and transfers by 
 \[\tau\!\uparrow^G_{H}:\mf{M}(H)\rtarr \mf{M}(G).\]
 \end{notn}

Mackey functors are required to satisfy the so-called double coset formula. Since our group $G=C_2\times C_2$ is abelian, this means that for any two $H_1$ and $H_2$ of the nontrivial cyclic subgroups, the restriction maps commute with the transfer maps, in the sense that
\begin{equation}\label{DoubleCoset}
r\!\downarrow^G_{H_1} \circ \ \tau\!\uparrow^G_{H_2} = \tau\!\uparrow_e^{H_1} \!\circ \  r\!\downarrow_e^{H_2}
\end{equation}

\begin{defn}
Given a surjection $\phi_N:G\rtarr G/N$ of groups with kernel $N \unlhd G$, there is a pullback for Mackey functors
\[ \phi_N^* : \mathrm{Mack}(G/N) \rtarr \mathrm{Mack}(G)\]
defined by 
\[ \phi_N^*(\mf{M})(H) = 
\begin{cases}
\mf{M} (H/N) & N\leq H \\
0 & N \not \leq H.
\end{cases}
\]
\end{defn}

In the Mackey functor literature, this pullback is known as inflation along the quotient $G\rtarr G/N$.

\begin{eg} Let 
\[ \MackC{\mf{M}(C_2)}{\mf{M}(e)}, \]
be a $C_2$-Mackey functor, where we  assume trivial Weyl group action for simplicity. Then, under the quotient map $\Kl \rtarr \Kl/R \iso C_2$, this pulls back to the $\Kl$-Mackey functor
\[
\xymatrix{
 & \mf{M}(R) \ar@/_1ex/[dr] & \\
 0 & 0 & \mf{M}(e)  \ar@/_1ex/[ul]  \\
  & 0. &
} \]
\end{eg}

\begin{notn}\label{InflNotn}
We will often encounter Mackey functors which are direct sums of inflations along the projections to different subgroups, and it will be convenient to use the notation
\[\phi^*_{LD}\mf{M} := \phi^*_L\mf{M} \oplus \phi_D^*\mf{M}, \qquad \phi^*_{LDR}\mf{M} := \phi^*_L\mf{M} \oplus \phi_D^*\mf{M} \oplus \phi_R^*\mf{M}. \]
\end{notn}


There is a related construction in the world of spectra. 
Given a surjection \mbox{$\phi:G \rtarr G/N$,}  there is a geometric pullback functor $\phi_N^*:\mathrm{Sp}^{G/N}\rtarr \mathrm{Sp}^G$ (\cite{LMS}*{Theorem~II.9.5}, \cite{H}*{Proposition~4.3}). For our purposes, the important property is its behavior on suspensions of Eilenberg-Mac~Lane spectra. This is given by
\begin{equation}
\phi_N^* (S^{V^G}\sma H_{G/N} \mf{M}) \simeq S^{V} \sma H_G \phi_N^* \mf{M}
\end{equation}
for $V\in RO(G)$ and $\mf{M}\in \mathrm{Mack}(G/N)$ (\cite[Proposition~4.2, Corollary~4.6]{H}).

\subsection{\for{toc}{Suspension, transfers, and restrictions}\except{toc}{Relationship between twisted (de)suspensions, transfers, and restrictions}}

Consider the $\Kl$-cofiber sequence $\Kl/R_+ \rtarr S^0 \rtarr S^\alpha$.
For any $\Kl$-spectrum $X$, this induces a cofiber sequence
\begin{equation}\label{AlSuspCofib} (\Kl/R_+\sma X)^\Kl \simeq X^R
  \xrtarr{\tau\uparrow_R^\Kl} X^\Kl \rtarr (\Sigma^\alpha X)^\Kl, \end{equation}
where the map $\tau$ is the transfer.
We similarly get that $(\Sigma^{-\alpha}X)^\Kl$ is the fiber of the resriction:
\begin{equation}\label{AlDesuspFib} (\Sigma^{-\alpha}X)^\Kl \rtarr X^\Kl \xrtarr{r\downarrow_R^\Kl} X^R.  \end{equation}

We have similar fiber sequences relating the $L$ transfer and restriction to (de)suspension by $\beta$, and the $D$ transfer and restriction to (de)suspension by $\gamma$.

\subsection{Anderson duality}

In this section, $G$ can be any finite group.
By Brown Representability in the category of $H\ulF$-modules,
the functor
\[ X \mapsto \Hom_{\F}(\pi_*^G X, \F)\]
on the category of $H\ulF$-modules is represented by some $H\ulF$-module, which we write $\F^{H\ulF}$. As in \cite[Lemma~3.1]{GM}, plugging in the $H\ulF$-modules $G/H_+\sma H\ulF$ shows that in fact $\F^{H\ulF}\simeq H\ulF^*$. The following more general result, whose proof was explained to us by John Greenlees, will be quite useful.

\begin{prop}\label{DualProp} Let $\ul{M}$ be an $\ulF$-module. Then
\[ \mpi_V(H\ul{M}^*) \iso \big( \mpi_{-V}H\ul{M}\big)^*.\]
\end{prop}

\begin{pf}
By Brown Representability in the category of $H\ulF$-modules,
the functor
\[ X \mapsto \Hom_{\F}(\pi_*^G (X\sma_{H\ulF}H\ul{M}), \F)\]
on the category of $H\ulF$-modules is represented by some $H\ulF$-module, which we write $\F^{H\ul{M}}$. By plugging in $X=G/H_+\sma H\ulF$, we see that $\F^{H\ul{M}} \simeq H\ul{M}^*$. In other words,
\[ [ X, H\ul{M}^*]^{H\ulF-\mathrm{mod}} \iso  \Hom_{\F}(\pi_*^G (X\sma_{H\ulF}H\ul{M}), \F).\]
Plugging in $X=S^V\sma H\ulF$ gives the result.
\end{pf}

\subsection{The slice filtration}\label{sec:slicefiltration}

We've already defined $X \geq n$ for a $G$-spectrum $X$ and we have a notion of ``less than''
as well.

\begin{defn} We say that $X< n$ if 
\[ [ S^{k\rho_H + r},X]^H =0 \]
for all $r\geq 0$ and all subgroups $H\leq G$ such that $k|H| \geq n$.
\end{defn}

In other words, $X<n$ if and only if the restriction $X\! \downarrow^G_{H}$ is less than $n$ for all proper subgroups $H<G$ and 
\[ [S^{k\rho_G+r},X]^G=0\]
for all $r\geq 0$ and $k\geq \frac{n}{|G|}$. More generally, restriction to subgroups is compatible with the slice dimension, in the following sense.

 \begin{prop}[\cite{H}*{Cor. 2.6}]\label{SliceRestriction} Suppose that $X\in \mathrm{Sp}^G$ satisfies $k \leq X\leq n$ and $H\leq G$. Then $k \leq X \! \downarrow^G_{H} \leq n$ as an $H$-spectrum.
\end{prop}

The following characterization of the subcategory $\tau_{\geq n}^G$ in
terms of connectivity of fixed points is useful.

\begin{thm}[\cite{HY}*{Corollary~2.9, Theorem~2.10}]\label{HYThm} Let $n\geq 0$. Then $X\geq n$ if and only if 
\[ \pi_k(X^H) = 0 \qquad \text{for\ } k < \frac{n}{|H|}.\]
\end{thm}

 An immediate corollary is

 \begin{cor}\label{nConnGEQn} If $n\geq 0$ and $X$ is $n$-connective, in the sense that $\pi_k(X^H)=0$ for all subgroups and all $k<n$, then $X\geq n$.
 \end{cor}

For a few values of $n$, the category of $n$-slices is well-understood.

\begin{prop}\label{CatSlices} \ 
\begin{enumerate}
\item \cite[Proposition~4.50]{Kervaire} $X$ is a $0$-slice if and only if $X\simeq H \mf{M}$ for $\mf{M}$ an arbitrary Mackey functor.
\item \cite[Proposition~4.50]{Kervaire} $X$ is a $1$-slice if and only if $X\simeq \SI^1 H \mf{M}$ for $\mf{M}$ a Mackey functor with injective restrictions.
\item \cite[Theorem 6-4]{U} $X$ is a $(-1)$-slice if and only if $X\simeq \SI^{-1} H \mf{M}$ for $\mf{M}$ a Mackey functor with surjective transfers.
\end{enumerate}
\end{prop}

Though these characterizations are not enough to determine all slices
in every case as the slice tower does not commute with taking ordinary
suspensions, it does commute with suspensions by the regular
representation of $G$.

\begin{prop}[\cite{Kervaire}*{Corollary 4.25}]\label{SlicesRho}
For any $k\in Z$,
\[ P^{k+|G|}_{k+|G|}(\SI^\rho X) \simeq \SI^\rho P^{k}_{k}(X).\]
\end{prop}

Additionally, we understand the relationship between the slice
filtration and taking pullbacks.

\begin{prop}[\cite{U}*{Corollary 4-5}]\label{SlicesPullbacks} Let $N\unlhd G$ be normal of index $k$ and let $X$ be a $G/N$-spectrum. Then
\[ \phi^*_N (P_n^n X) \simeq P_{kn}^{kn} ( \phi_N^* X).\]
In particular, the pullback of an $n$-slice is a $kn$-slice.
\end{prop}

 \begin{prop}\label{dSliceEquiv} Let $d\in \Z$ and let 
 \[ X \xrtarr{f} Y \rtarr Z\]
 be a fiber sequence of $G$-spectra such that $P^d_d(Z)\simeq \ast \simeq P^d_d(\SI^{-1}Z)$. Then $f$ induces an equivalence on $n$-slices. 
 \end{prop}
 
 \begin{pf}
 This follows from \cite[Proposition~2.32]{W}.
 \end{pf}

\subsection{Review of Holler-Kriz}

In \cite{HK}, the authors compute the homotopy of $(\SI^V H\ulF)^G$ for any elementary abelian group $G$. Their answer is given as the Poincar\'{e} series of the graded $\F$-vector space.

\begin{thm}[\cite{HK}*{Section 6}]\label{HKtwo} Let $ \ell, n \geq 0$ and $i, j \geq 1$. 
The Poincar\'{e} series for $\pi_*((\SI^V H\ulF)^\Kl)$ is
\begin{enumerate}
\item $V=0$: \quad $1$
\item $V=n\al$: \quad $1+x+\dots+x^n$ 
\item $V=-j\al$: \quad  $x^{-j} + \dots + x^{-3}+ x^{-2}$
\item $V=n\al+\ell\be$: \quad $(1+\dots + x^n) \cdot (1+\dots+x^\ell)$
\item $V=n\al-j\be$ : \quad $(1+\dots+x^n) \cdot (x^{-j}+\dots + x^{-2})$
\item $V=-i\al-j\be$ : \quad $(x^{-i}+\dots+x^{-2}) \cdot (x^{-j}+\dots + x^{-2})$
\end{enumerate}
\end{thm}

If either $i$ or $j$ is equal to $1$, then the above series should be interpreted as zero.
The answer is more complicated when all three nontrivial irreducible representations are involved, so we state those cases separately.

\begin{thm}[\cite{HK}*{Section 6}]\label{HKPos} Let $  \ell ,  m ,   n \geq 1$. 
The Poincar\'{e} series for $\pi_*((\SI^{\ell\al+m\be+n\ga} H\ulF)^\Kl)$ is
\[ (1+\dots+x^\ell)(1+\dots+x^m) + x(1+\dots+x^{\ell+m})(1+\dots+x^{n-1})\]
\end{thm}

Expanded out, this polynomial can be described as follows, assuming $\ell\leq m\leq n$: The constant coefficient is $1$. Then the coefficients increase by 2 until $x^\ell$. Thereafter, they increase by 1 until $x^m$. They then stay constant until $x^n$, and finally decrease (by 1) to 1, which is the coefficient of $x^{\ell+m+n}$.

\begin{thm}[\cite{HK}*{Section 6}] Let $  \ell , m \geq 1 $. If $ k\geq 2$,  
then the Poincar\'{e} series for $\pi_*((\SI^{\ell\al+m\be-k\ga} H\ulF)^\Kl)$ is
\[ \left( \frac1{x^k} + \dots +\frac1x\right)(1+x+\dots+x^{k-2}) + x^k(1+\dots+x^{\ell-k})(1+\dots+x^{m-k})
\]
In the case $k=1$, the series is
\[ x(1+\dots+x^{\ell-1})(1+\dots+x^{m-1}).\]
\end{thm}

\begin{thm}[\cite{HK}*{Section 6}]\label{HKljk} Let $ j, k, \ell \geq 1$. 
Then the Poincar\'{e} series for $\pi_*((\SI^{\ell\al-j\be-k\ga} H\ulF)^\Kl)$ is
\begin{enumerate}
\item
\[ \frac1{x^{j+k-\ell}}(1+\dots+x^{j-\ell-2})(1+\dots+x^{k-\ell-2}) + \frac1{x^{\ell+1}}(1+\dots+x^\ell)(1+\dots+x^{\ell-1})
\]
if $j,k \geq \ell+1$ or
\item
\[ \frac1{x^j}(1+\dots+x^{j-2})(1+\dots+x^{\ell-k}) + \frac1{x^k}(1+\dots+x^{\ell-1})(1+\dots+x^{k-1})
\]
if $\ell\geq k$.
\end{enumerate}
\end{thm}

Swapping the role of $j$ and $k$ gives the case $\ell\geq j$ in \autoref{HKljk}.

\begin{thm}[\cite{HK}*{Section 6}]\label{HKijk} Let $i, j,  k\geq 1$. 
Then the  Poincar\'{e} series for $\pi_*((\SI^{-i\al-j\be-k\ga} H\ulF)^\Kl)$ is
\[ \frac1{x^{i+j+k}} \Big[ (1+x+\dots+x^{j+k-2})(1+\dots+x^{i-2}) + x^{i-1}(1+\dots+x^{k-1})(1+\dots+x^{j-1}) \Big]
\]
\end{thm}

\begin{cor}\label{HKk3} Let $ k\geq 1$. 
The Poincar\'{e} series for $\pi_*((\SI^{k-k\rho} H\ulF)^\Kl)$ is
\[ \frac1{x^{3k}} \Big[ (1+x+\dots+x^{2k-2})(1+\dots+x^{k-2}) + x^{k-1}(1+\dots+x^{k-1})^2 \Big]
\]
\end{cor}

\section{Review of $G=C_2$}\label{sec:C2}

A Mackey functor for the group $C_2$ may be depicted by the Lewis diagram
\[ \MackC{M(G)}{M(e)}, \]
where we have omitted the $C_2$-action on $M(e)$.

\subsection{The main players}\label{sec:C2Mackey}

\begin{eg} The constant Mackey functor is $\ulF=\raisebox{2em}{\MackCAr{\F}{\F}{1}{0}}$.
\end{eg}

\begin{eg} The geometric Mackey functor is $\underline{g}=\phi_{C_2}^*(\F)=\raisebox{2em}{\MackC{\F}{0}}$. Since $\underline{g}(e)=0$, it follows that $(H\underline{g})^e\simeq \ast$. Smashing the cofiber sequence
\[ (C_2)_+ \rtarr S^0 \rtarr S^\sigma\] 
with $H\underline{g}$ 
 implies that 
\[\Sigma^{k\sigma} H\underline{g} \simeq H\underline{g} \qquad \text{and} \qquad \Sigma^{k\rho}H\underline{g} \simeq \Sigma^k H\underline{g}  . \]
Thus, using either \autoref{SlicesPullbacks} or \autoref{SlicesRho}, it follows that $\SI^k H\ulg$ is a $2k$-slice.
\end{eg}

\begin{eg} The free Mackey functor is $\underline{f}=\raisebox{2em}{\MackC{0}{\F}}$. This is relevant because 
\begin{equation}\label{ulfModel}
\Sigma^{1-\sigma}H\ulF \simeq H\underline{f}.
\end{equation}
\end{eg} 

Note that these Mackey functors sit in an exact sequence
\[ \underline{f} \into \ulF \onto \underline{g}.\]
The resulting cofiber sequence 
\begin{equation} \label{Fcofiber}
 H\underline{f} \rtarr H\ulF \rtarr H\underline{g}
 \end{equation}
can be used to compute the homotopy of $\SI^{k\rho} H \ulF$.

\begin{prop}\label{HtpykrhoC2F} For $k\geq 0$, the nontrivial homotopy Mackey functors of $\SI^{k\rho}H_{C_2} \ulF$ are
\[ \mpi_i\left( \SI^{k\rho} H_{C_2} \ulF\right) \iso 
\begin{cases}
\ulF & i = 2k \\
\ulg & i \in [k,2k-1]
\end{cases} \]
\end{prop}

\begin{pf}
This follows by induction from repeated use of the cofiber sequence
\[ \SI^{(j-1)\rho + 2} H\ulF \simeq  \SI^{j \rho} H\ulf \rtarr \SI^{j \rho} H \ulF \rtarr \SI^{j \rho} H\ulg \simeq \SI^j H\ulg,\]
where $j\geq 1$.
\end{pf}

\begin{eg} The opposite to the constant Mackey functor is $\ulF^*=\raisebox{2em}{\MackCAr{\F}{\F}{0}{1}}$. We again have a twisting
\begin{equation}\label{ulfopModel}
\Sigma^{1-\sigma}H\ulf \simeq H\ulF^*.
\end{equation}
\end{eg}

We also have the  exact sequence of Mackey functors
\[ \ulg \into \ulF^* \onto \ulf.\]
The resulting cofiber sequence 
\begin{equation} \label{Fopcofiber}
 H\ulg \rtarr H\ulF^* \rtarr H\ulf
 \end{equation}
can be used to compute the homotopy of $\SI^{-k\rho} H\ulF^*$ (which also follows from \autoref{HtpykrhoC2F} by \autoref{DualProp}).

\begin{prop}\label{HtpykrhoC2F*} For $k\geq 0$, the nontrivial homotopy Mackey functors of $\SI^{-k\rho}H_{C_2} \ulF^*$ are
\[ \mpi_{-i}\left( \SI^{-k\rho} H_{C_2} \ulF^*\right) \iso 
\begin{cases}
\ulF^* & i = 2k \\
\ulg & i \in [k,2k-1]
\end{cases} \]
\end{prop}

Together, \autoref{HtpykrhoC2F} and \autoref{HtpykrhoC2F*} combine to give the $RO(C_2)$-graded homotopy Mackey functors of $H\ulF$, which we display in \autoref{fig:C2}.

\begin{figure}[h]
\caption{
\label{fig:C2}
$\mpi_{n+k\sigma}H_{C_2}\ulF$ }
\begin{center}
\begin{tikzpicture}
\draw [step = 1] (-5,-3) grid (4,6);
\foreach \x in {-5,...,3}
 \node[anchor=north] at (\x+0.5,-3) {\x};
\foreach \y in {-3,...,5}
 \node[anchor=east] at (-5,\y+0.5) {\y};
\draw [ultra thick,<->] (0,-3) -- (0,6);
\draw [ultra thick,<->] (-5,0) -- (4,0);
\node at (-5.8,5.5) {$\sigma$};
\draw [thick,->] (-5.6,5.2) -- (-5.6,5.8);
\foreach \z in {0,...,3}
  \node at (0.5+\z,0.5-\z) {$\ulF$};
\foreach \z in {-1,...,-3}
  \node at (0.5,0.5+\z) {$\ulg$};
\foreach \z in {-2,-3}
  \node at (1.5,0.5+\z) {$\ulg$};
\node at (2.5,-2.5) {$\ulg$};
 \node at (-0.5,1.5) {$\ulf$};
\foreach \z in {2,...,5}
  \node at (0.5-\z,0.5+\z) {$\ulF^*$};
\foreach \z in {3,...,5}
  \node at (-1.5,0.5+\z) {$\ulg$};
\foreach \z in {4,5}
  \node at (-2.5,0.5+\z) {$\ulg$};
\node at (-3.5,5.5) {$\ulg$};
\end{tikzpicture}
\end{center}
\end{figure}

\subsection{\for{toc}{The slice tower for $\Sigma^n H F_2$}\except{toc}{The slice tower for $\Sigma^n H_{C_2} \ulF$}}

\begin{eg} The spectrum $H\ulF$ is a 0-slice by \autoref{CatSlices}. 
\end{eg}

\begin{eg} The spectrum $\Sigma^1H\ulF$ is a $1$-slice by \autoref{CatSlices}.
\end{eg}

\begin{eg} The spectrum $\Sigma^2 H\ulF$ is a $2$-slice, since
\[ \Sigma^2 H\ulF \simeq \Sigma^\rho \left( \Sigma^{1-\sigma} H\ulF\right) \simeq \Sigma^\rho H\underline{f}.\]
\end{eg}

\begin{eg} The spectrum $\Sigma^3 H\ulF$ is a $3$-slice, since
\[ \Sigma^3 H\ulF \simeq \Sigma^\rho \left( \Sigma^{2-\sigma} H\ulF\right) \simeq \Sigma^\rho \Sigma^1H\underline{f}.\]
Since $\Sigma^1 H\underline{f}$ is a 1-slice, the claim follows. 
\end{eg}

\begin{eg} The spectrum $\Sigma^4 H\ulF$ is a $4$-slice, since
\[ \Sigma^4 H\ulF \simeq \Sigma^\rho \left( \Sigma^{3-\sigma} H\ulF\right) \simeq \Sigma^\rho \Sigma^2H\underline{f}
\simeq \Sigma^{2\rho} \Sigma^{1-\sigma} H\underline{f} \simeq \Sigma^{2\rho} H\ulF^*.\]
Since $ H\ulF^*$ is a 0-slice, the claim follows. 
\end{eg}

\begin{eg} The spectrum $\Sigma^5 H\ulF$ has both a $5$-slice and a $6$-slice, since
\[ \Sigma^5 H\ulF \simeq \Sigma^\rho \left( \Sigma^{4-\sigma} H\ulF\right) \simeq \Sigma^\rho \Sigma^3H\underline{f}
\simeq \Sigma^{2\rho} \Sigma^{2-\sigma} H\underline{f} \simeq \Sigma^{2\rho} \Sigma^1 H\ulF^*.\]
Now the slice tower for $\Sigma^1 H\ulF^*$ is the suspension of \eqref{Fopcofiber}:
\[  \Sigma^1 H\underline{g} \rtarr \Sigma^1 H\ulF^* \rtarr \Sigma^1 H\underline{f}.\]
The left spectrum is a 2-slice, while the right one is a 1-slice. Thus the slice tower for $\SI^5 H\ulF$ is 
\[ P^6_6 = \SI^3 H\ulg \rtarr \SI^5 H\ulF \rtarr \SI^{2\rho+1}H\ulF = P_5^5.\]
\end{eg}

More generally, we have

\begin{thm}\label{SlicesSigmanHF}
The slice tower of $\Sigma^n H\ulF$, for $n\geq 4$, is
\medskip

\hspace{-1em}
\begin{minipage}{0.4\textwidth}
\centerline{$n$ even} \vspace{-1em}
\[\scalebox{0.8}
{ \xymatrix @C=-4em{ 
P^{2n-4}_{2n-4} = \Sigma^{n-2} H\underline{g} \ar[r] & \Sigma^n H\ulF \ar[d]\\
P^{2n-6}_{2n-6} = \Sigma^{n-3} H\underline{g} \ar[r] & \Sigma^{(n-1)+\sigma} H\ulF \ar[d]\\
		& \ar@{.}[d] \\
P^{n+2}_{n+2} = \Sigma^{\frac{n}{2}+1} H\underline{g} \ar[r] & \Sigma^{\left(\frac{n}{2} -3\right) \rho + 6} H\ulF \ar[d]\\
		& \mbox{ $P^n_n = \Sigma^{(\frac{n}2-2) \rho + 4}H\ulF \simeq \Sigma^{\frac{n}{2}\rho } H{\ulF^*}  $\hspace{10em}}  }
		}\]
\end{minipage} \hspace{0.12\textwidth}
\begin{minipage}{0.45\textwidth}
\centerline{$n$ odd} \vspace{-1em}
\[\scalebox{0.8}{ \xymatrix @C=-6em{ 
P^{2n-4}_{2n-4} = \Sigma^{n-2} H\underline{g} \ar[r] & \Sigma^n H\ulF \ar[d]\\
P^{2n-6}_{2n-6} = \Sigma^{n-3} H\underline{g} \ar[r] & \Sigma^{(n-1)+\sigma} H\ulF \ar[d]\\
		& \ar@{.}[d] \\
P^{n+1}_{n+1} = \Sigma^{\frac{n+1}{2}} H\underline{g} \ar[r] & \Sigma^{\left(\frac{n+1}{2} -3\right) \rho + 5} H\ulF \ar[d]\\
		& \mbox{$P^n_n = \Sigma^{\frac{n-3}2\rho+3}H\ulF\simeq \Sigma^{\frac{n-1}{2}\rho + 1} H\mf{f} $\hspace{10em}} }
		}\]
\end{minipage}
\end{thm}

\begin{pf}

The $2\rho$-suspension of $H\ulg \rtarr H\ulF^* \rtarr H\ulf$ is
\[ \SI^2 H\ulg \rtarr \SI^4 H\ulF \rtarr \SI^{\rho+2}H\ulF.\]
The theorem is obtained by repeated application of suspensions of this cofiber sequence.
\end{pf}

\begin{cor}\label{SlicesSigmanHf} The $C_2$-spectrum $\Sigma^n H\ulf$ is an $n$-slice for $n=0,1,2$. If $n\geq 3$, then $ n \leq \Sigma^n H\ulf \leq 2n-2$ and $P^{2n-2}_{2n-2}(\Sigma^n H\ulf) \simeq \Sigma^{n-1}H\ulg$. 
\end{cor}

\begin{pf} This follows from the fiber sequence
\[ \Sigma^{n-1} H\ulg \rtarr \Sigma^n H\ulf \rtarr \Sigma^n H\ulF\]
and \autoref{SlicesSigmanHF}. Indeed, the slice tower is given by augmenting the slice tower for $\Sigma^n H\ulF$ with the above fiber sequence.
\end{pf}

\section{Mackey functors for $\Kl=C_2\times C_2$}\label{sec:KMackey}

A Lewis diagram for a Mackey functor over the Klein-four group takes the shape

\[ \MackK{M(\Kl)}{M(L)}{M(D)}{M(R)}{M(e)} \]
We have not drawn in the $C_2$-actions on the intermediate groups or the $\Kl$-action on $M(e)$ (these actions are trivial in all of our examples).
In the examples below, we only draw restriction or transfer maps that are nonzero.

\begin{eg}
We have the constant Mackey functor
\[ \ulF = \raisebox{3.8em}{\xymatrix{
 & \F \ar[dl]_{1} \ar[d]_{1} \ar[dr]^1 & \\
 \F \ar[dr]_1& \F \ar[d]_1  & \F \ar[dl]^1,   \\ 
  & \F  &
}} \]
as well as its dual
\[ \ulF^* = 
\raisebox{3.8em}{\xymatrix{
 & \F \ & \\
 \F \ar[ur]^1& \F \ar[u]^1  & \F \ar[ul]_1,   \\ 
  & \F  \ar[ul]^1  \ar[u]^1 \ar[ur]_1 &
}} 
\]
\end{eg}

\begin{prop}\label{RhoDesusp} $\Sigma^{-\rho}H\ulF \simeq \Sigma^{-4} H\ulF^*$.
\end{prop}

\begin{pf}
Restricting to $L$, say, we have 
\[ \iota_L^* (\Sigma^{-\rho}H_\Kl \ulF) \simeq \Sigma^{-2-2\sigma}H_{C_2}\ulF \simeq \Sigma^{-4}\Sigma^{2-2\sigma}H_{C_2}\ulF \simeq \Sigma^{-4} H_{C_2}\ulF^*.\]
 The same argument applies to the restriction to $D$ and $R$. 
\autoref{HKijk} gives that $(\Sigma^{-\rho}H\ulF)^\Kl \simeq \Sigma^{-4} H\F$. 

The transfer map from $R$ to $\Kl$ fits into a fiber sequence 
\[ (\SI^{-2\rho}H_R\ulF)^R \simeq (\Kl/R_+ \sma \SI^{-\rho}H_\Kl \ulF)^\Kl \rtarr (\SI^{-\rho}H_\Kl \ulF)^\Kl \rtarr (\SI^{\al-\rho}H_\Kl \ulF)^\Kl.\]
By \autoref{HKtwo} this becomes a fiber sequence
\[ \SI^{-4} H\F \rtarr \SI^{-4}H\F \rtarr \ast,\]
so that the transfer map is an equivalence. By symmetry, we find that the other transfer maps are equivalences as well.

Similarly, the restriction from $\Kl$ to $R$ fits into a fiber sequence
\[ (\SI^{-\rho-\al}H_\Kl \ulF)^\Kl \rtarr (\SI^{-\rho}H_\Kl \ulF)^\Kl \rtarr (\Kl/R_+ \sma \SI^{-\rho} H_\Kl \ulF)^\Kl \simeq (\SI^{-2\rho}H_R\ulF)^R.\]
By \autoref{HKijk}, the spectrum $(\SI^{-\rho-\al}H_\Kl \ulF)^\Kl$ has $\pi_{-4}\iso \F \iso \pi_{-5}$. It follows from the long exact sequence in homotopy that the restriction map must be zero.
\end{pf}

\begin{eg} The geometric Mackey functor is
\[ \ulg := \phi_\Kl^*(\F) = 
 \raisebox{4em}{\xymatrix{
 & \F  & \\
 0 & 0  & 0,   \\ 
  & 0  &
}} \]
\end{eg}

\begin{eg} The free Mackey functor is
\[ \ulf :=  
 \raisebox{4em}{\xymatrix{
 & 0 & \\
 0 & 0  & 0,   \\ 
  & \F  &
}} \]
\end{eg}

Unlike the case for $G=C_2$, the $\Kl$-spectrum $H_\Kl \ulf$ is {\em not} equivalent to $\SI^V H_\Kl \ulF$ for any $V$.

\begin{eg}\label{mDefn}
\[ \mf{m} := 
\raisebox{3.8em}{\xymatrix{
 & \F \ar[dr]^1 \ar[dl]_1 \ar[d]_1 & \\
 \F   & \F  & \F  .  \\
  & 0 &
}}
\]
and
\[  \mf{m}^* :=
\raisebox{3.8em}{\xymatrix{
 & \F  & \\
 \F  \ar[ur]^1 & \F \ar[u]^1  & \F \ar[ul]_1  .  \\
  & 0 &
}}
\]
\end{eg}

\begin{eg}\label{wDefn}
\[  \mf{w} :=
\raisebox{3.8em}{\xymatrix{
 & 0  & \\
 \F  \ar[dr]^1 & \F \ar[d]^1  & \F \ar[dl]_1  .  \\
  & \F &
}}
\]
and
\[ \mf{w}^* := 
\raisebox{3.8em}{\xymatrix{
 & 0 & \\
 \F   & \F  & \F  .  \\
  &  \F \ar[ur]^1 \ar[ul]_1 \ar[u]_1  &
}}
\]
\end{eg}

\begin{eg}\label{mgDefn}
\[ \mf{mg} := 
\raisebox{3.8em}{\xymatrix{
 & \F \oplus \F \ar[dr]^{p_2} \ar[dl]_{p_1} \ar[d]_{\nabla} & \\
 \F   & \F  & \F  .  \\
  & 0 &
}}
\]
and
\[  \mf{mg}^* :=
\raisebox{3.8em}{\xymatrix{
 & \F \oplus \F  & \\
 \F  \ar[ur]^{i_1} & \F \ar[u]^\Delta  & \F \ar[ul]_{i_2}  .  \\
  & 0 &
}}
\]
\end{eg}

\begin{prop}\label{mProp}
There are equivalences
\begin{enumerate}
\item
$ \SI^{-\rho} H\mf{m} \simeq \SI^{-2} H \mf{mg}^*$ 
\item
 $ \SI^{\rho} H\mf{m}^* \simeq \SI^{2} H \mf{mg}$
\end{enumerate}
\end{prop}

\begin{pf}
We prove the second statement. The first follows in a similar way, or by citing \autoref{DualProp}.
Consider the (nonsplit) short exact sequence
\[ \phi^*_L \ulF^* \into \mf{m}^* \onto \phi_{DR}^* \ulf.\]
This gives a nonsplit cofiber sequence 
\[ \SI^2 H\phi_L^*\ulf \simeq \SI^{\rho} H\phi^*_L \ulF^* \rtarr \SI^{\rho} H\mf{m}^* \rtarr \SI^{\rho} H\phi_{DR}^* \ulf \simeq \SI^2 H\phi_{DR}^* \ulF.\]
It follows that $\SI^\rho H\mf{m}^*\simeq \SI^2 H\mf{E}$, for some nonsplit extension $\mf{E}$ of $\phi_{DR}^*\ulF$ by $\phi^*_L \ulf$. 
But we can similarly express $\mf{E}$ as an extension of $\phi_{LR}^*\ulF$ by $\phi^*_D \ulf$.
The only possibility is $\mf{E}\iso \mf{mg}$.
\end{pf}

\begin{eg}\label{WDefn}
\[ \mf{W} := 
\raisebox{3.8em}{\xymatrix{
 & \F \oplus \F \oplus \F & \\
 \F  \ar[ur]^{i_1} \ar[dr]_1 & \F \ar[u]^{i_2} \ar[d]_1  & \F \ar[ul]_{i_3}  \ar[dl]^1.  \\
  &  \F   &
}}
\]
and
\[  \mf{W}^* :=
\raisebox{3.8em}{\xymatrix{
 & \F\oplus \F\oplus \F \ar[dl]_{p_1} \ar[d]_{p_2} \ar[dr]^{p_3}  & \\
 \F  & \F   & \F  .  \\
  & \F  \ar[ul]^1 \ar[u]_1 \ar[ur]_1 &
}}
\]
\end{eg}

\section{\for{toc}{The slices of $\SI^n H_\Kl F_2$}\except{toc}{The slices of $\SI^n H_\Kl \ulF$}}\label{sec:slices}

\begin{prop}\label{SliceBounds} For $n\geq 0$ and $\Kl=C_2\times C_2$, the $\Kl$-spectrum $\SI^n H\ulF$ satisfies $n\leq \SI^n H\ulF \leq 4(n-3)$.
\end{prop}

\begin{pf} The lower bound follows by \autoref{nConnGEQn}. For the upper bound, note that after restricting to the trivial subgroup, the spectrum is an $n$-slice. 
By \autoref{SlicesSigmanHF}, the restriction to a cyclic subgroup is bounded above by $2n-4$ if $n\geq 4$ (it is an $n$-slice if $n=0,\dots,3$). It therefore remains to check that 
\[  \pi_{k\rho_\Kl+r}(\SI^n H\ulF) =    [ S^{k\rho_\Kl+r},\SI^n H\ulF]^\Kl=0\]
for $r\geq 0$ and $4k > 4(n-3)$. In other words, the homotopy groups $\pi_r^\Kl \big(\SI^{n-k\rho_\Kl}H\ulF\big)$ must vanish if $r\geq 0$ and $k>n-3$.
This follows from \autoref{HKk3}.
\end{pf}

Moreover, we know a priori in which dimensions the slices appear.

\begin{prop}\label{EvenSliceBounds} For $n\geq 0$ and $\Kl=C_2\times C_2$, all slices of the $\Kl$-spectrum $\SI^n H\ulF$ above level $n$ are even slices. Furthermore, if $n\geq 4$, then all slices above level $2n-4$ occur only in dimensions that are multiples of $4$.
\end{prop}

\begin{pf}
Since the restriction to any cyclic subgroup  is bounded above by $2n-4$ by \autoref{SlicesSigmanHF},
it follows that 
all slices of $\Sigma^n H\ulF$ above
dimension $2n-4$ must be geometric. Thus, we know further that the
only nontrivial slices of $\Sigma^n H\ulF $ above dimension $2n-4$ are
$4k$-slices. 

Similarly, by \autoref{SlicesSigmanHF}, the restriction of $\SI^n H\ulF$ to any cyclic subgroup has only even slices, except for possibly the $n$-slice. Thus any odd slices above level $n$ must be geometric. But geometric $\Kl$-spectra only have slices in dimension a multiple of $|\Kl|=4$. It follow that all odd slices above level $n$ must be trivial.
\end{pf}

\subsection{The $n$-slice}

We will use the following recursion to establish the bottom slices of $\SI^n H\ulF$.

\begin{prop}\label{RecursionProp}
Let $n\geq 7$. Then
\[ P_k^k ( \SI^n H\ulF) \simeq \SI^\rho P_{k-4}^{k-4}(\SI^{n-4} H\ulF)\]
for $k\in [ n, 2n-7]$.
\end{prop}

\begin{pf}
By \autoref{SlicesRho}, we have
\[ P_k^k(\SI^n H\ulF) \simeq \SI^\rho P^{k-4}_{k-4}(\SI^{n-\rho}H\ulF) \simeq \SI^\rho P^{k-4}_{k-4}(\SI^{n-4}H\ulF^*).\]
Thus it suffices to compare the $(k-4)$-slices of $\SI^{n-4}H\ulF^*$ and $\SI^{n-4}H\ulF$.

The short exact sequences of Mackey functors
\[ 0 \rtarr \mf{m}^* \rtarr \ulF^* \rtarr \ulf \rtarr 0 \]
and 
\[ 0 \rtarr \ulf \rtarr \ulF \rtarr \mf{m} \rtarr 0\]
give the following diagram of fiber sequences
\[\xymatrix{
\SI^j H\mf{m}^* \ar[r] \ar[d] & \SI^j H\ulF^* \ar@{=}[d] \ar[r] & \SI^j H\ulf \ar[d] \\
\fib(\lambda) \ar[r] \ar[d] & \SI^j H\ulF^* \ar[r]^{\lambda} \ar[d] & \SI^j H\ulF \ar[d] \\
\SI^{j-1} H\mf{m} \ar[r] & \ast \ar[r] & \SI^j H\mf{m}
}\]
Then $\fib(\lambda)$ is $j-1$-connective, and the underlying spectrum of $\fib(\lambda)$ is contractible. By \autoref{HYThm}, it follows that $\fib(\lambda)\geq 2j-2$ as long as $j\geq 1$. Similarly, $\SI \fib(\lambda) \geq 2j$. By \cite[Corollary~4.17]{HHR},
if tollows that 
$\lambda$ induces an equivalences of slices below $2j-2$. Taking $j=n-4$ gives the result.
\end{pf}

Note that the above argument, using only the fiber sequence $\SI^j H\mf{m}^* \rtarr \SI^j H\ulF^* \rtarr \SI^j H\ulf$, gives the following result that will be employed below.

\begin{prop}\label{HalfRecursionProp}
Let $n\geq 7$. Then
\[ P_k^k ( \SI^n H\ulF) \simeq \SI^\rho P_{k-4}^{k-4}(\SI^{n-4} H\ulf)\]
for $k\in [ n, 2n-5]$.
\end{prop}

\begin{prop}\label{nSliceProp} For $n\geq 0$ and $\Kl=C_2\times C_2$, the bottom slice of $\SI^n H\ulF$ is 
\[ P_n^n (\SI^n H\ulF) \simeq 
\begin{cases}
\SI^n H\ulF & n\in [0,4] \\
\SI^{\frac{n}4\!\rho}H\ulF^* & n \equiv 0 \pmod4, n\geq 4 \\
\SI^{\frac{n-1}4\!\rho+1}H\ulf & n \equiv 1 \pmod4, n\geq 4 \\
\SI^{\frac{n-2}4\!\rho+2}H\ulf & n \equiv 2 \pmod4, n\geq 4 \\
\SI^{\frac{n-3}4\!\rho+3}H\ulF & n \equiv 3 \pmod4, n\geq 4
\end{cases}
\]
\end{prop}

\begin{pf}
By \autoref{RecursionProp}, it suffices to establish the base cases, in which $n\leq 7$. These are given in \autoref{sec:towers} below.
\end{pf}

\subsection{The $4k$-slices}

\begin{thm}\label{4kSlicesThm} 
For all $n > 4$,
\[P^{4k}_{4k}(\Sigma^n H\ulF) = \begin{cases} \Sigma^k H\ulg^{2(n-k)-5} &
4k \in [2n-3,4n-12]  \\
\Sigma^{k\rho} H \Big(\phi_{LDR}^*\,(\ulF^*) \oplus 
\mf{g}^{4k-n-2}\Big) & 4k \in [n+2,2n-4] \\
\Sigma^{k\rho} H \mf{mg}^* & 4k = n +1  \\
\Sigma^{k\rho} H \ulF^* & 4k = n  \\
* & \mbox{otherwise} 
\end{cases}\]
\end{thm}

\begin{pf} 
The above formula agrees with 
\autoref{SlicesSigmanHF} upon restriction to the cyclic subgroups. 
To determine the $\Kl$-fixed points, we use that
\[ P^{4k}_{4k}(\SI^n H\ulF) \simeq \SI^{k\rho} H \mf{\pi}_0 \SI^{n-k\rho}H\ulF\]
by repeated application of \autoref{SlicesRho}. The fixed points are then given by \autoref{4kSliceFxdPtsHelper}.
\medskip

It remains to consider the restriction and transfer maps if $4k \in [n,2n-4]$. The restriction maps to the subgroup $R$, for instance, 
 fit into fiber sequences
\[ \big( \SI^{n-k\rho-\al} H \ulF\big)^\Kl \rtarr \big( \SI^{n-k\rho} H \ulF\big)^\Kl \xrtarr{res} \big( \SI^{n-k\rho} H \ulF\big)^R.\]
Fixing $k>1$, we argue by induction on $n$ that \autoref{4kSliceRestrictionHelper} implies that the long exact sequence in homotopy splits into a series of short exact sequences of $\F$-vector spaces
\[ 0 \to \mf{\pi}_{k\rho-n+1}^R H\ulF \into  \mf{\pi}_{k\rho+\al-n}^\Kl H\ulF \onto \mf{\pi}_{k\rho-n}^\Kl H\ulF \to 0,\]
linked together by the null restriction map. Since $4k\in [n,2n-4]$, it follows that $2k+2\leq n \leq 4k$.

The base case for our induction argument is the case $n=2k+2$, so that $4k$ is $2n-4$. In this base case, the left term $\mf{\pi}_{k\rho-2k-1}^R H\ulF$ vanishes, and the other two terms are both of dimension $2k-1$. This establishes the base case.

For the induction step, we suppose that $\mf{\pi}_{k\rho-n+1}^R H\ulF \into  \mf{\pi}_{k\rho+\al-n}^\Kl H\ulF$ is injective. It follows that we have an exact sequence
\[ 0 \rtarr \F \into \F^{4k-n+2} \to \F^{4k-n+1},\]
which shows that the map on the right must be surjective.

A similar argument shows that the transfer map from the subgroup $R$, say,  up to $\Kl$ is injective. We argue by (downward) induction on $n$
that \autoref{4kSliceTransferHelper} implies that the long exact sequence in homotopy for the fiber sequence
\[  \big( \SI^{n-k\rho} H \ulF\big)^R \xrtarr{tr}
 \big( \SI^{n-k\rho} H \ulF\big)^\Kl 
\rtarr \big( \SI^{n-k\rho+\al} H \ulF\big)^\Kl \]
splits into a series of short exact sequences
\[ 0 \to \mf{\pi}_{k\rho-n}^R H\ulF \into  \mf{\pi}_{k\rho-n}^\Kl H\ulF \onto \mf{\pi}_{k\rho-\al-n}^\Kl H\ulF \to 0.\]
The base case is $n=2k+1$, so that $4k=2n-2$. In this case,
$\mf{\pi}_{k\rho-2k-1}^RH\ulF=0$, and the other two terms are both of
the same dimension (using \autoref{4kSliceTransferHelper}).

For the induction step, we suppose that the transfer map $\mf{\pi}_{k\rho-n}^R H\ulF \into  \mf{\pi}_{k\rho-n}^\Kl H\ulF$ is injective. It follows that we have an exact sequence
\[ 0 \rtarr \F \into \F^{4k-n+1} \to \F^{4k-n},\]
which shows that the map on the right must be surjective.

It remains to show that the transfer maps are linearly independent if $4k \geq n+2$ and have distinct images if $4k=n+1$. 
In the case $4k=n+1$, consider the exact sequence of Mackey functors
\[ \mf{\pi}_0\big( \Kl/R_+ \sma \SI^{4k-1-k\rho} H\ulF \big) \rtarr \mf{\pi}_0\big(  \SI^{4k-1-k\rho} H\ulF \big) \rtarr \mf{\pi}_0\big(  \SI^{4k-1-k\rho+\al} H\ulF \big)\]
The left Mackey functor vanishes at $L$ and $D$, and the $R$-transfer map in the middle Mackey functor is in the image of the left Mackey functor. 
Thus to see that the $L$ or $D$ transfers in the middle Mackey functor have image distinct
 from that of the $R$ transfer, it suffices to show that the $L$ or $D$ transfer in the right Mackey functor is nontrivial.
But a similar argument to that for the transfer maps above shows that the Mackey functor on the right is $\ulF^*$, so we are done. By symmetry, we similarly conclude that the images of the $L$ and $D$ transfers are distinct. 

Finally, if $4k\geq n+2$, to show additionally that the three transfers are linearly independent, we consider the exact sequence
\[ \pi_0\big( \Kl/D_+ \sma \SI^{n-k\rho+\al+\be} H\ulF \big) \rtarr {\pi}_0\big(  \SI^{n-k\rho+\al+\be} H\ulF \big) \rtarr {\pi}_0\big(  \SI^{n-k\rho+\al+\be+\ga} H\ulF \big)\]
This is an exact sequence of the form
\[ \F \rtarr \F^{4k-n-1} \rtarr \F^{4k-n-2},\]
so that the first map cannot be zero. We conclude that the $D$ transfer map is nonzero after factoring out the $L$ and $R$ transfers, so that the three transfer maps are independent if $4k\geq n+2$.
\end{pf}

The following lemmas are direct consequences of \autoref{HKijk}.

\begin{lem}\label{4kSliceFxdPtsHelper}
For $n\geq 4$, we have
\[ \pi_0 \big( \SI^{n-k\rho} H \ulF\big)^\Kl \iso \begin{cases} \F^{2(n-k)-5} & 4k \in [2n-3,4n-12] \\ 
\F^{4k-n+1} & 4k \in [n,2n-4] \\
0 & \text{else} \end{cases} \]
\end{lem}

\begin{pf}
The dimension of the fixed points is given by the coefficient of $x^{k-n}$
in the Poincar\'e series of \autoref{HKk3}. Equivalently, the dimension is given by
 the coefficient of $x^{4k-n}$ in the polynomial
 \[ p(x) = (1+x+\dots+x^{2k-2})(1+\dots+x^{k-2}) + x^{k-1}(1+\dots+x^{k-1})^2\]
 This polynomial can be described as follows:
 The constant coefficient is $1$. Then the coefficients increase by 1 until $(2k-1)x^{2k-2}$  and then decrease by 2 until $1\cdot x^{3k-3}$. In other words, the coefficient of $x^i$ is
 \[ \begin{cases} 
 i+1 & 0\leq i \leq 2k-2 \\
 6k-2i-5 & 2k-1 \leq i \leq 3k-3 \\
 0 & \text{else}.
 \end{cases} \]
 Plugging in $i=4k-n$ gives the result.
\end{pf}

\begin{lem}\label{4kSliceRestrictionHelper}
For $n\geq 4$,
we have
\[ \pi_0 \big( \SI^{n-k\rho-\al} H \ulF\big)^\Kl \iso \begin{cases}
  \F^{2(n-k)-5} & 4k \in [2n-4,4n-12] \\ 
\F^{4k-n+2} & 4k \in [n-1,2n-6] \\
0 & \text{else} \end{cases} \]
\end{lem}

\begin{pf}
The dimension of the fixed points here is given by the coefficient of $x^{k-n}$
in the Poincar\'e series of \autoref{HKijk}. Equivalently, the dimension is given by
 the coefficient of $x^{4k-n+1}$ in the polynomial
 \[ p(x) = (1+x+\dots+x^{2k-2})(1+\dots+x^{k-1}) + x^{k-1}(1+\dots+x^{k-1})^2\]
The polynomial is nearly the same as that from
\autoref{4kSliceFxdPtsHelper} and can be described as follows:
 The constant coefficient is $1$. Then the coefficients increase by 1
 until $(2k-1)x^{2k-2}$, remain constant for the term
 $(2k-1)x^{2k-1}$, and  then decrease by 2 until $1\cdot x^{3k-2}$. In other words, the coefficient of $x^i$ is
 \[ \begin{cases} 
 i+1 & 0\leq i \leq 2k-2 \\
 6k-2i-3 & 2k-1 \leq i \leq 3k-2 \\
 0 & \text{else}.
 \end{cases} \]
 Plugging in $i=4k-n+1$ gives the result.
\end{pf}

\begin{lem}\label{4kSliceTransferHelper}
For $n>4$, 
we have
\[ \pi_0 \big( \SI^{n-k\rho+\al} H \ulF\big)^\Kl \iso \begin{cases} \F^{2(n-k)-5} & 4k \in [2n-2,4n-12] \\ 
\F^{4k-n} & 4k \in [n+1, 2n-4] \\
 0 & \text{else} \end{cases} \]
\end{lem}

\begin{pf}
Since $n > 4$, $k > 1$ so the dimension of the fixed points in this
case is still given by the coefficient of $x^{k-n}$
in the Poincar\'e series of \autoref{HKijk}. Equivalently, the dimension is given by
 the coefficient of $x^{4k-n-1}$ in the polynomial
 \[ p(x) = (1+x+\dots+x^{2k-2})(1+\dots+x^{k-3}) + x^{k-2}(1+\dots+x^{k-1})^2\]
The polynomial is similar to those in 
\autoref{4kSliceFxdPtsHelper} and
\autoref{4kSliceRestrictionHelper} and can be described as follows:
 The constant coefficient is $1$. Then the coefficients increase by 1
 until $(2k-2)x^{2k-3}$, decrease by 1 for the single term
 $(2k-3)x^{2k-2}$, and  then decrease by 2 until $1\cdot x^{3k-4}$. In other words, the coefficient of $x^i$ is
 \[ \begin{cases} 
 i+1 & 0\leq i \leq 2k-3 \\
 6k-2i-7 & 2k-2 \leq i \leq 3k-4 \\
 0 & \text{else}.
 \end{cases} \]
 Plugging in $i=4k-n-1$ gives the result.
\end{pf}

\subsection{The $(4k+2)$-slices}\label{Sec:4k+2Slices}

In this section, we obtain the $4k+2$-slices. We begin with the top such slices.

\begin{prop} Let $n\geq 8$ be even. Then 
\[ P^{2n-6}_{2n-6}(\SI^n H\ulF) \simeq \SI^{(\frac{n}2-2)\rho+1}H\phi_{LDR}^* \ulf.\]
\end{prop}

\begin{pf}
By \autoref{HalfRecursionProp}, we have
\[ P^{2n-6}_{2n-6}(\SI^n H\ulF) \simeq \SI^\rho P^{2n-10}_{2n-10}(\SI^{n-4}H\ulf).\]
The short exact sequence
\[ 0 \rtarr \ulf \rtarr \ulF \rtarr \mf{m} \rtarr 0\]
gives rise to the fiber sequence 
\[ \SI^{n-5} H \mf{m} \rtarr  \SI^{n-4} H\ulf \rtarr \SI^{n-4} H\ulF.\]
By \autoref{EvenSliceBounds}, we have
\[ P^{2n-10}_{2n-10}(\SI^{n-5} H\ulF) \simeq \ast \simeq P^{2n-10}_{2n-10}(\SI^{n-4} H\ulF).\]
It then follows from \autoref{dSliceEquiv} that $\SI^{n-5}H\mf{m}\rtarr \SI^{n-4}H\ulf$ induces an isomorphism on $(2n-10)$-slices.

The short exact sequence 
\[ 0 \rtarr \mf{m} \rtarr \phi^*_{LDR}\ulF \rtarr \ulg^2 \rtarr 0 \]
gives a fiber sequence 
\[ \SI^{n-6} H \ulg^2 \rtarr \SI^{n-5} H\mf{m} \rtarr \SI^{n-5} H \phi^*_{LDR}\ulF.\]
If $n\geq 8$ is even, then by \autoref{SlicesSigmanHF}, we get that 
\[P^{2n-10}_{2n-10}(\SI^{n-5}H\mf{m}) \simeq \SI^{(\frac{n}2-3)\rho+1}H\phi^*_{LDR}\ulf.\]
\end{pf}

\begin{prop} Let $n \geq 5$ be odd. Then 
\[ P^{2n-4}_{2n-4}(\SI^n H\ulF) \simeq \SI^{(\frac{n+1}2-2)\rho+1}H\phi_{LDR}^* \ulf.\]
\end{prop}

\begin{pf}
For $n=5$, this is given in \autoref{sigma5tower} below. 
We have
\[ P^{2n-4}_{2n-4}(\SI^n H\ulF) \simeq P^{2n-4}_{2n-4}(\SI^{\rho+n-4}H\ulF^*) \simeq \SI^\rho P^{2n-8}_{2n-8}(\SI^{n-4}H\ulF^*).\]
The short exact sequence
\[ 0 \rtarr \ulg \rtarr \ulF^* \rtarr \mf{w}^* \rtarr 0\]
gives a fiber sequence
\[ \SI^{n-4}H\ulg \rtarr \SI^{n-4}H\ulF^* \rtarr \SI^{n-4} H\mf{w}^*.\]
It follows that, if $n\geq 5$, then $\SI^{n-4}H\ulF^* \rtarr \SI^{n-4}H\mf{w}^*$ induces an equivalence on $(2n-8)$-slices.

Next, the short exact sequence
\[ 0 \rtarr \mf{w}^* \rtarr \mf{W}^* \rtarr \ulg^3 \rtarr 0\]
yields a fiber sequence
\[ \SI^{n-5}H\ulg^3 \rtarr \SI^{n-4}H\mf{w}^* \rtarr \SI^{n-4}H\mf{W}^*.\]
It follows that, if $n\geq 7$, then $\SI^{n-4}H\mf{w}^* \rtarr \SI^{n-4}H\mf{W}^*$ induces an equivalence on $(2n-8)$-slices.

Finally, consider the short exact sequnce
\[ 0 \rtarr \phi^*_{LDR}\ulF \rtarr \mf{W}^* \rtarr \ulf \rtarr 0.\]
This gives a fiber sequence
\[ \SI^{n-4}H\phi^*_{LDR}\ulF \rtarr \SI^{n-4} H \mf{W}^* \rtarr \SI^{n-4}H\ulf.\]
The restriction of $\SI^{n-4}H\ulf$ to either $L$, $D$, or $R$ has no slices above level $2n-12$, and similarly for $\SI^{n-5}H\ulf$. Thus the $(2n-8)$-slice of $\SI^{n-4}H\ulf$ (and $\SI^{n-5}H\ulf$) must be geometric. Since $2n-8$ is not a multiple of 4, we conclude that the $(2n-8)$-slices are trivial. By \autoref{dSliceEquiv}, we conclude that $ \SI^{n-4}H\phi^*_{LDR}\ulF \rtarr \SI^{n-4} H \mf{W}^*$ induces an equivalence on $(2n-8)$-slices.
We are now done by \autoref{SlicesSigmanHF}.
\end{pf}

\begin{prop}\label{4k2SlicesProp} Let $4k+2 \in (n,2n-4]$. Then
\[ P^{4k+2}_{4k+2}(\SI^n H\ulF) \simeq \SI^{k\rho+1} H\phi_{LDR}^*\ulf.\]
\end{prop}

\begin{pf}
If $4k+2\leq 2n-7$, this follows from \autoref{RecursionProp} and the base cases discussed in \autoref{sec:towers}. This leaves only the cases of $4k+2=2n-6$ if $n$ is even, or $4k+2=2n-4$ if $n$ is odd. 
These cases are handled in the two preceding propositions.
\end{pf}

\section{\for{toc}{The slice towers of $\SI^n H_\Kl F_2$}\except{toc}{The slice towers of $\SI^n H_\Kl \ulF$}}
\label{sec:towers}

We determine the slice towers of 
$\Sigma^n H\ulF$ 
for $0 \leq n\leq 8$.

\begin{eg} 
$H\ulF$ is a zero-slice. 
\end{eg}

\begin{eg}
$\Sigma^1 H\ulF$ is a 1-slice since the restriction maps are injective.
\end{eg}

\begin{eg}
$\Sigma^2 H\ulF$ is a 2-slice. Since this is true upon restriction to each of the proper subgroups, it suffices to show that 
\[ [S^{n\rho+r},\Sigma^2 H\ulF]=0 \]
for $n>0$ and $r\geq 0$. 
This follows from \autoref{HKijk}. 
\end{eg}

\begin{eg}
$\Sigma^3 H\ulF$ is a 3-slice. Since this is true upon restriction to each of the proper subgroups, it suffices to show that 
\[ [S^{n\rho+r},\Sigma^3 H\ulF]=0 \]
for $n,r>0$. This follows from \autoref{HKijk}. Alternatively, $\Sigma^{-\rho}\Sigma^3H\ulF\simeq \SI^{-1} H\ulF^*$ is a $(-1)$-slice according to \autoref{CatSlices} since the transfer maps are surjective. It follows that $\SI^3 H\ulF$ is a 3-slice.

\end{eg}

\begin{eg} By \autoref{RhoDesusp}, $\Sigma^4 H\ulF \simeq \Sigma^\rho H\ulF^*$. It follows that $\SI^4 H\ulF$ is a 4-slice.
\end{eg}

\begin{eg}\label{sigma5tower}
Consider the short exact sequences of Mackey functors
\[ 0\rtarr \ulg \rtarr \ulF^* \rtarr \mf{w}^* \rtarr 0\]
and
\[ 0 \rtarr \phi^*_{LDR}\ulf \rtarr \mf{w}^* \rtarr \ulf \rtarr 0,\]
 where $\mf{w}^*$ is defined in \autoref{wDefn}.
The resulting cofiber sequences produce the
slice tower for $\Sigma^5 H\ulF \simeq \Sigma^\rho \Sigma^1 H\ulF^*$:
\[ \xymatrix{
P^8_8 = \Sigma^{2} H\mf{g} \ar[r] & \SI^{\rho+1} H\ulF^* \simeq \SI^5 H\ulF \ar[d] \\
P^6_6 = \Sigma^{\rho+1}H \phi_{LDR}^*\ulf   \ar[r] & 
\Sigma^{\rho+1}H \mf{w}^* \ar[d]   \\
  & P_5^5 = \Sigma^{\rho+1} H\ulf,
 }\]
\end{eg}

\begin{eg}\label{sigma6tower}
Suspending the slice tower for $\Sigma^5 H\ulF$ gives the tower for $\SI^6 H\ulF$:
\[ \xymatrix{
P^{12}_{12} = \Sigma^{3} H\mf{g} \ar[r] & \Sigma^{\rho+2} H\ulF^* \simeq \SI^6 H\ulF \ar[d] \\
P^8_8 = \Sigma^{\rho+2}H \phi_{LDR}^*\ulf   \ar[r] &
\Sigma^{\rho+2}H \mf{w}^* \ar[d]   \\
  & P_6^6 = \Sigma^{\rho+2} H\ulf.
 }\]

\begin{lem} $\SI^2 H\ulf$ is a 2-slice.
\end{lem}

\begin{pf}
The short exact sequence of Mackey functors
\[ 0 \rtarr \ulf \rtarr \ulF \rtarr \mf{m} \rtarr 0,\]
where $\mf{m}$ is defined in \autoref{mDefn},
gives a cofiber sequence
\[ \SI^1 H\mf{m} \rtarr \SI^2 H\ulf \rtarr \SI^2 H\ulF.\]
The spectrum $\SI^2 H\ulF$ is a 2-slice, and the cofiber sequence
\[ \SI^1H \phi_{L,D}^* \ulf \rtarr \SI^1 H\mf{m} \rtarr \SI^1 H \phi_R^* \ulF\]
shows that $\SI^1 H\mf{m}$ is also a 2-slice.
\end{pf}
\end{eg}

\begin{eg} For $\SI^7 H\ulF \simeq \SI^\rho\SI^3 H\ulF^*$, we have fiber sequences as in
\[ \xymatrix{
P^{16}_{16} = \Sigma^{4} H\mf{g} \ar[r] & \SI^{\rho+3}H\ulF^* \simeq \SI^7 H\ulF \ar[d] \\
 P_{12}^{12} = \SI^{3} H \ulg^3 \ar[r] & \Sigma^{\rho+3}H \mf{w}^*  \ar[d]  \\
P_{10}^{10}= \SI^{\rho+3} H  \phi_{LDR}^*\, \ulF  \ar[r] & \SI^{\rho+3} H\mf{W}^* \ar[d] \\
 P_8^8 = \SI^{\rho+2} H\mf{m} \ar[r] & P_7^8 = \Sigma^{\rho+3} H\ulf \ar[d] \\
  & P_7^7 = \SI^{\rho+3} H\ulF,
 }\]
 where $\mf{W}^*$ is defined in \autoref{WDefn}.
 Note that $\SI^2 H\mf{m}$ is a 4-slice, as 
 $\SI^2 H\mf{m} \simeq \SI^\rho H\mf{mg}^*$ according to \autoref{mProp}.
\end{eg}

\begin{eg} For $\SI^8 H\ulF \simeq \SI^\rho\SI^4 H\ulF^*$, we have fiber sequences as in
\[ \xymatrix{
P^{20}_{20} = \Sigma^{5} H\mf{g} \ar[r]  & \SI^{\rho+4}H\ulF^* \simeq \SI^8 H\ulF \ar[d] \\
P_{16}^{16} = \SI^{4} H \ulg^3 \ar[r]  
& \Sigma^{\rho+4}H \mf{w}^*  \ar[d]  \\
P_{12}^{12} =   \SI^{3\rho} H \big( \phi_{LDR}^*\, \ulF^*\oplus \ulg^2\big) \ar[r] & \SI^{\rho+4}H\mf{W}^* \ar[d]  \\
 & P^{10} \SI^8 H\ulF
 }\]
The twelve slice is given by \autoref{4kSlicesThm}. 

\begin{prop} The slice section $P^{10}\SI^8 H\ulF$ is equivalent to $\SI^{\rho+3}C$, where $C$ is the cofiber of $H\ulF \rtarr H\phi^*_{LDR}\,\ulF$.
\end{prop}

\begin{pf}
The cofiber $C$ has homotopy Mackey functors $\mf{\pi}_1(C) \iso \ulf$ and $\mf{\pi}_0(C) \iso \mf{g}^2$. Thus the $\rho$-suspension of the Postnikov sequence is a fiber sequence
\[ \SI^{\rho+2}H\mf{g}^2 \rtarr \SI^{\rho+4}H\ulf \rtarr \SI^{\rho+3} C.\]
On the other hand, the short exact sequence $\phi^*_{LDR}\,\ulF \into \mf{W}^* \rtarr \ulf$ gives a fiber sequence
\[ \SI^{3\rho} H\phi^*_{LDR}\, \ulF^* \simeq \SI^{\rho+4}H\phi^*_{LDR}\ulF \rtarr \SI^{\rho+4} H\mf{W}^* \rtarr \SI^{\rho+4}H \ulf.\]
Since the 12-slice is the sum of the left terms in these two sequences and $\SI^{\rho+4} H\mf{W}^*$ is $P^{12}\SI^8 H \ulF$, the octahedral axiom gives a cofiber sequence 
\[ P_{12}^{12}(\SI^8 H\ulF) \rtarr P^{12} \SI^8 H\ulF \rtarr \SI^{\rho+3}C.\]
But $8 \leq \SI^{\rho+3}C \leq 10$, so we are done.
\end{pf}

\begin{rem}
The $\Kl$-spectrum $\SI^{\rho+3}C$ is the first example of a $\Kl$-spectrum that is not an $RO(\Kl)$-graded suspension of an Eilenberg-Mac~Lane spectrum and yet which occurs as a slice or slice section in the tower for $\SI^n H\ulF$. Indeed, the restriction of $C$ to each cyclic subgroup is $\SI^1H\ulf$. But if $\SI^V H\mf{M}$ restricts to $\SI^1 H\mf{M}$ for each cyclic subgroup, it must be that $V=1$. As $C$ has a nontrivial $\mf{\pi}_0$, it cannot be of the form $\SI^V H\mf{M}$.
\end{rem}

Thus the slice tower is given by
\[ \xymatrix{
P^{20}_{20} = \Sigma^{5} H\mf{g} \ar[r]  & \SI^{\rho+4}H\ulF^* \simeq \SI^8 H\ulF \ar[d] \\
P_{16}^{16} = \SI^{4} H \ulg^3 \ar[r]  
& \Sigma^{\rho+4}H \mf{w}^*  \ar[d]  \\
P_{12}^{12} =   \SI^{3\rho} H \big( \phi_{LDR}^*\, \ulF^*\oplus \ulg^2\big) \ar[r] & \SI^{\rho+4}H\mf{W}^* \ar[d]  \\
  P_{10}^{10} = \SI^{\rho+3} H \phi_{LDR}^*\, \ulF \ar[r] & \SI^{\rho+3} C \ar[d] \\
  & P_8^8 = \SI^{\rho+4} H\ulF,
 }\]
\end{eg}

For the higher suspensions, we do not know the slice tower explicitly. We give a diagram of fibrations which is close to the slice tower in the next two examples.

\begin{eg}\label{SI9tower} For $\SI^9 H\ulF \simeq \SI^\rho\SI^5 H\ulF^*$, we have fiber sequences as in
\[ \xymatrix{
P^{24}_{24} = \Sigma^{6} H\mf{g} \ar[rr] & & \SI^{\rho+5}H\ulF^* \simeq \SI^9 H\ulF \ar[d] \\
P_{20}^{20} = \SI^{5} H \ulg^3 \ar[rr] &  
& \Sigma^{\rho+5}H \mf{w}^*  \ar[d]  \\
\fbox{A}\  \SI^4 H\ulg^3 \ar[r] &  \SI^{\rho+5} H  \phi_{LDR}^*\, \ulF  \ar[d] \ar[r] & \SI^{\rho+5}H\mf{W}^* \ar[dd]  \\
 & P_{14}^{14} = \SI^{3\rho+1} H  \phi_{LDR}^*\, \ulf  &  \\
\fbox{B}\  \SI^{\rho+3} H\ulg^2 \ar[rr] &  
& \Sigma^{\rho+5} H\ulf \ar[d] \\
 &  \fbox{C}\ \SI^{\rho+4} H \phi_{LDR}^*\, \ulF   \ar[r] & \SI^{\rho+4} C \ar[d] \\
  &\fbox{D}\  \SI^{\rho+2} \ulg \ar[r] &  \SI^{\rho+5} H\ulF \ar[d] \\
   & P^{10}_{10} =  \SI^{2\rho+1}H\phi_{LDR}^* \ulf \ar[r] & \SI^{2\rho+1} H\mf{w}^* \ar[d] \\
    & & P^{9}_{9} = \SI^{2\rho+1} H\ulf,
 }\]
This is not quite the slice tower. According to \autoref{4kSlicesThm}, the 16-slice is the sum $\fbox{A}\vee\fbox{B}$. Similarly, the 12-slice is the sum $\fbox{C}\vee \fbox{D}$.
\end{eg}

\begin{eg}\label{SI10tower} For $\SI^{10} H\ulF \simeq \SI^\rho\SI^6 H\ulF^*$, we have fiber sequences as in
\[ \xymatrix{
P^{28}_{28} = \Sigma^{7} H\mf{g} \ar[rr] & & \SI^{\rho+6}H\ulF^* \simeq \SI^{10} H\ulF \ar[d] \\
P_{24}^{24} = \SI^{6} H \ulg^3 \ar[rr] &  
& \Sigma^{\rho+6}H \mf{w}^*  \ar[d]  \\
\fbox{A}\  \SI^5 H\ulg^3 \ar[r] &  \SI^{\rho+6} H  \phi_{LDR}^*\, \ulF  \ar[d] \ar[r] & \SI^{\rho+6}H\mf{W}^* \ar[dd]  \\
 & \fbox{C}\ \SI^{3\rho+2} H  \phi_{LDR}^*\, \ulf  &  \\
 \fbox{B}\ \SI^{\rho+4} H\ulg^2 \ar[rr] &  
& \Sigma^{\rho+6} H\ulf \ar[d] \\
\fbox{D}\ \SI^{3\rho+1} H\ulg^3 \ar[r] & \SI^{\rho+5} H \phi_{LDR}^*\, \ulF   \ar[r] \ar[d] & \SI^{\rho+5} C \ar[dd] \\
 & P^{14}_{14} = \SI^{3\rho+1} H \phi_{LDR}^* \ulf & \\
  &\fbox{E}\ \SI^{\rho+3} \ulg \ar[r] &  \SI^{\rho+6} H\ulF \ar[d] \\
   & P^{12}_{12} =  \SI^{2\rho+2}H\phi_{LDR}^* \ulf \ar[r] & \SI^{2\rho+2} H\mf{w}^* \ar[d] \\
    & & P^{10}_{10} = \SI^{2\rho+2} H\ulf,
 }\]
According to \autoref{4kSlicesThm}, the 20-slice is the sum $\fbox{A}\vee\fbox{B}$ and the 16-slice is the sum $\fbox{C}\vee \fbox{D}\vee \fbox{E}$.
\end{eg}

\section{Homotopy Mackey functor computations}\label{sec:Homotopy}

Here we collect some computations of homotopy Mackey functors of various twisted Eilenberg-Mac~Lane spectra.

\begin{thm}\label{Htpy-krhoF*} 
For all $k\geq 1$, the nontrivial homotopy Mackey functors of $\SI^{-k\rho}H\ulF^*$ are
\[\mpi_{-i}(\Sigma^{-k\rho} H\ulF^*) = \left\{\begin{array}{cl}
\ulF^* & i = 4k \\
\mf{mg}^* & i = 4k-1 \\
\phi^*_{LDR}\ulF^* \oplus \ulg^{4k-2-i} & i \in [2k , 4k-2] \\
\ulg^{2(i-k)+1} & i \in [k,2k-1]. \\
\end{array}\right.\]
\end{thm}

\begin{pf} In the proof of \autoref{4kSlicesThm}, we computed the
  following homotopy group Mackey functors:

\[\mpi_{-n}(\Sigma^{-k\rho} H\ulF) = \left\{\begin{array}{cl}
\ulF^* & n = 4k \\
\mf{mg}^* & n = 4k-1 \\
\phi^*_{LDR}\ulF^* \oplus \ulg^{4k-2-n} & 4k \in [n+2,2n-4] \\
\ulg^{2(n-k)-5} & 4k \in [2n-3, 4n-12]. 
\end{array}\right.\]

We know that 
\[\mpi_{-i}(\SI^{-k\rho}H\ulF^*) = \mpi_{-(i+4)}(\SI^{-(k+1)\rho}
H\ulF).\]
The result follows by setting $n=i+4$ and replacing $k$ by $k+1$.
\end{pf}

\autoref{DualProp} then gives the following result, which gives the homotopy Mackey functors of the bottom slice in the case $n\equiv 0,3\pmod4$.

\begin{cor}\label{HtpykrhoF} 
For all $k\geq 1$, the nontrivial homotopy Mackey functors of $\SI^{k\rho}H\ulF$ are
\[\mpi_i(\Sigma^{k\rho} H\ulF) = \left\{\begin{array}{cl}
\ulF & i = 4k \\
\mf{mg} & i = 4k-1 \\
\phi^*_{LDR}\ulF \oplus \ulg^{4k-2-i} & i \in [2k , 4k-2] \\
\ulg^{2(i-k)+1} & i \in [k,2k-1]. \\
\end{array}\right.\]
\end{cor}

For the homotopy of the bottom slice in the remaining two cases, namely $n\equiv 1,2\pmod4$, we need some auxiliary computations.

\begin{prop}\label{Htpykrhom} Let $k\geq 1$. The nontrivial homotopy Mackey functors of $\SI^{k\rho} H\mf{m}$ are
\begin{enumerate}
\item $\mf{\pi}_{2k} (\Sigma^{k\rho}\mf{m}) \iso \phi^*_{LDR}\ulF$ 
\item $\mf{\pi}_{i} (\Sigma^{k\rho}\mf{m}) \iso \ulg^3$ for $i\in [k+1,2k-1]$ 
\item $\mf{\pi}_{k} (\Sigma^{k\rho}\mf{m}) \iso \ulg$ 
\end{enumerate}
\end{prop}

\begin{pf} 
The short exact sequence
\[ \phi^*_{LDR}\ulf \into \mf{m} \onto \ulg \]
gives a cofiber sequence 
\[ \SI^{(k-1)\rho+2}H\phi^*_{LDR}\ulF \simeq \SI^{k\rho} H \phi^*_{LDR} \ulf \rtarr \SI^{k\rho} H \mf{m} \rtarr \SI^{k\rho} H \ulg \simeq \SI^k H\ulg.\]
The result now follows from \autoref{HtpykrhoC2F}.
\end{pf}

The same argument, using instead the short exact sequence $\phi^*_{LDR}\ulf \into \mf{m} \onto \ulg^2$, applies to show

\begin{prop}\label{Htpykrhomg} Let $k\geq 1$. The nontrivial homotopy Mackey functors of $\SI^{k\rho} H\mf{mg}$ are
\begin{enumerate}
\item $\mf{\pi}_{2k} (\Sigma^{k\rho}\mf{mg}) \iso \phi^*_{LDR}\ulF$ 
\item $\mf{\pi}_{i} (\Sigma^{k\rho}\mf{mg}) \iso \ulg^3$ for $i\in [k+1,2k-1]$ 
\item $\mf{\pi}_{k} (\Sigma^{k\rho}\mf{mg}) \iso \ulg^2$ 
\end{enumerate}
\end{prop}

The homotopy of the bottom slice when $n \equiv 1,2\pmod 4$ is now given by the following result.

\begin{cor}\label{Htpykrhof}
For all $k\geq 1$, the nontrivial homotopy Mackey functors of $\SI^{k\rho}H\ulf$ are
\[\mpi_i(\Sigma^{k\rho} H\ulf) = \left\{\begin{array}{cl}
\ulF & i = 4k \\
\mf{mg} & i = 4k-1 \\
\phi^*_{LDR}\ulF \oplus \ulg^{4k-2-i} & i \in [2k+1 , 4k-2] \\
\ulg^{2(i-k-1)} & i \in [k+2,2k]. \\
\end{array}\right.\]
\end{cor}

\begin{pf}
We have the cofiber sequence
\[ \SI^{k\rho} H \ulf \rtarr \SI^{k\rho} H\ulF \rtarr \SI^{k\rho} H
\mf{m}\]
arising from the short exact sequence of Mackey functors. 
From the long exact sequence in homotopy we get the desired homotopy
for $i \in [2k+1, 4k]$ because $\mf{\pi_i}(\SI^{k\rho}H\mf{m}) = 0$ for
$i > 2k$ by \autoref{Htpykrhom}.

When $i = 2k$, in the long exact sequence we have
\[ \mf{\pi_{2k}}(\SI^{k\rho}H\ulf) \hookrightarrow \phi^*_{LDR}\ulF
  \oplus \ulg^{2k-2} \rightarrow \phi^*_{LDR}\ulF.\]
On subgroups of size $2$, the map on the right is an isomorphism and
thus maps $\phi^*_{LDR}\ulF$ isomorphically to the target forcing
$\mf{\pi_{2k}}(\SI^{k\rho}H\ulf) = \ulg^{2k-2}$.

For $ i \in [k+2, 2k-1]$ we have
\[\mf{\pi_i}(\SI^{k\rho}H\ulf) \rightarrow \ulg^{1+2(i-k)} \rightarrow
\ulg^3 \]
and thus $\mf{\pi_i}(\SI^{k\rho}H\ulf) = \ulg^{j}$, where $j \geq
2(i-k-1)$. 

To show that $j \leq 2(i-k-1)$, we use the cofiber sequence
\[ \SI^{(k-1)\rho+4}H\ulF \simeq \SI^{k\rho} H \ulF^* \rtarr
\SI^{k\rho} H \ulf \rtarr \SI^{k\rho+1} H\mf{m}^*\simeq
\SI^{(k-1)\rho+3} H\mf{mg}.\]
For $i \in [k+3, 2k-1]$ we have the following in the long exact
sequence in homotopy:
\[ \ulg^{2(i-k)-5} \rightarrow \mf{\pi_i}(\SI^{k\rho} H\ulf) \rightarrow
  \ulg^3\]
and thus, $\mf{\pi_i}(\SI^{k\rho}H\ulf) = \ulg^{j}$ where $j \leq
2(i-k-1)$. When $i = k+2$, we have $\mf{\pi_i}(\SI^{k\rho}H\ulf) =
\ulg^2$ as desired since $\mf{\pi_i}(\SI^{(k-1)\rho + 4}H\ulF) = 0$ 
and $\mf{\pi_i}(\SI^{(k-1)\rho+3}H\mf{mg}) = \ulg^2$ by
\autoref{Htpykrhomg}. For $i<k+2$ we can see from either long exact
sequence that $\mf{\pi_i}(\SI^{k\rho}H\ulf) = 0$.
\end{pf}

The homotopy of the slices in dimension congruent to 2 modulo 4 is much simpler.

\begin{prop}
The nontrivial homotopy Mackey functors of $\SI^{k\rho +1}
H\phi^*_{LDR}\ulf$ are
\[\mpi_i(\Sigma^{k\rho +1} H\phi^*_{LDR}\ulf) =
\left\{ \begin{array}{cl}
\phi^*_{LDR}\ulF & i = 2k+1 \\
\ulg^3 & i \in [k+2, 2k] \\
\end{array}\right.\]
\end{prop}

\begin{pf}
This follows directly from \autoref{HtpykrhoC2F}, given that $\SI^{\rho_{C_2}} H_{C_2}\ulf\simeq \SI^2 H_{C_2}\ulF$.
\end{pf}

\section{The slice spectral sequence}\label{sec:sliceSS}

The Mackey functor-valued slice spectral sequence for $\SI^n H\ulF$ must recover that the only nontrivial homotopy Mackey functor is $\mpi_n(\SI^n H\ulF) \iso \ulF$. This Mackey functor already occurs in the bottom slice, and all higher slices get wiped out by the spectral sequence. To a large extent, the answer forces many of the differentials. Furthermore, the slice spectral sequence must restrict to recover the slice spectral sequence on each cyclic subgroup, which further allows us to deduce many differentials. In practice, only a few differentials require further argument. We discuss this in several examples.

\begin{eg}
The first example which has a nontrivial slice tower is $\SI^5 H\ulF$. 
The slices are
\[ P_5^5(\SI^5 H\ulF) = \SI^{\rho+1} H\ulf, \qquad P_6^6(\SI^5 H\ulF) \simeq \SI^3 H\phi^*_{LDR}\ulF, \qquad P_8^8(\SI^5 H\ulF) = \SI^2 H\ulg.\]
Thus the $6$ and $8$-slices are Eilenberg-Mac Lane. The homotopy Mackey functors of the $5$-slice are given in \autoref{Htpykrhof}

In the slice spectral sequence, the only possibility is that we have differentials
\[ d_1:\mf{\pi}_4(P_5^5)  \iso \mf{mg} \into \mf{\pi}_3(P_6^6) = \phi_{LDR}^*\, \ulF .\]
and
\[ d_2:\mf{\pi}_3(P^6_6)/d_1 \xrtarr{\iso} \mf{\pi}_2(P^8_8) \iso \ulg.\]
\end{eg}

The slice spectral sequence for $\SI^6 H\ulF$ is just the suspension of that for $\SI^5 H\ulF$. Those for $\SI^7 H\ulF$ and $\SI^8 H\ulF$ are not much more complicated. There is only one possible pattern of differentials, which is displayed in \autoref{78Fig}.

\begin{eg}
In the slice spectral sequence for $\SI^9 H\ulF$, displayed in \autoref{910Fig}, almost all differentials are forced by the fact that only $\pi_9(P_9^9 \SI^9 H\ulF)\iso \ulF$ can survive the spectral sequence. The sole exception is that the summand $\ulg$ of 
\[\mf{\pi}_6(P_9^9 \SI^9 H\ulF) \iso \phi_{LDR}^* \ulF \oplus \ulg\]
 can support either a $d_3$ to $\mf{\pi}_5(P_{12}^{12}\SI^9 H\ulF) \iso \ulg^3$ or a $d_5$ to $\mf{\pi}_5(P_{14}^{14}\SI^9 H\ulF) \iso \ulg^3$.

To see that it must in fact support the shorter $d_3$, we use that the map 
\[\SI^9 H\ulF \simeq \SI^{\rho+5}H\ulF^* \rtarr \SI^{\rho+5}\ulf\] 
(see \autoref{SI9tower}) induces an equivalence on $9$, $10$, and $12$-silces. Since $\SI^{\rho+5}H\ulf$ only has nontrivial $\mf{\pi}_9$ and $\mf{\pi}_8$ by \autoref{Htpykrhof}, there must be a $d_3$ in the slice spectral sequence for $\SI^{\rho+5}H\ulf$ in order to wipe out the $\mf{\pi}_6$ and $\mpi_5$. 
\end{eg}

\begin{eg}
Most differentials in the slice spectral sequence for $\SI^{10}H\ulF$, which is displayed in \autoref{910Fig}, are forced by the fact that only $\mpi_{10}(P_{10}^{10}\SI^{10}H\ulF)$ survives.

To see that $d_2:\mpi_6(P_{10}^{10}) \rtarr \mpi_5(P_{12}^{12})$ is injective, we use that the map 
\[\SI^{10}H\ulF \simeq \SI^{\rho+6} H\ulF^* \rtarr \SI^{2\rho+2}H\mf{w}^*\]
 (see \autoref{SI10tower}) induces an equivalence on $10$ and $12$ slices. The cofiber sequence 
\[ \SI^{\rho+4}H\ulF \rtarr \SI^{2\rho}H\mf{w}^* \rtarr \SI^{2\rho+1}H\ulg\]
shows that $\mpi_6(\SI^{2\rho+2}H\mf{w}^*)=0$, which forces the claimed $d_2$-differential.

Similarly, the $\ulg$ summand of $\mpi_7(P_{10}^{10})$ supports a $d_4$ to the $14$-slice. This
can be seen by using the map to $\SI^{\rho+5}C$, which induces an equivalence of slices up to level $14$. The cofiber sequence
\[ \SI^{\rho+4}H\ulf \rtarr \SI^{\rho+3}C \rtarr \SI^{\rho+3}\ulg^2\]
shows that $\mpi_5(\SI^{\rho+3}C)=0$ and $\mpi_6(\SI^{\rho+3}C)\iso \ulg^2$. This forces the claimed $d_4$-differential.
\end{eg}

\newpage

\begin{figure}
\caption{\label{78Fig}The slice spectral sequence over $C_2$ and $C_2\times C_2$, $n=7,8$}
\end{figure}
\includegraphics{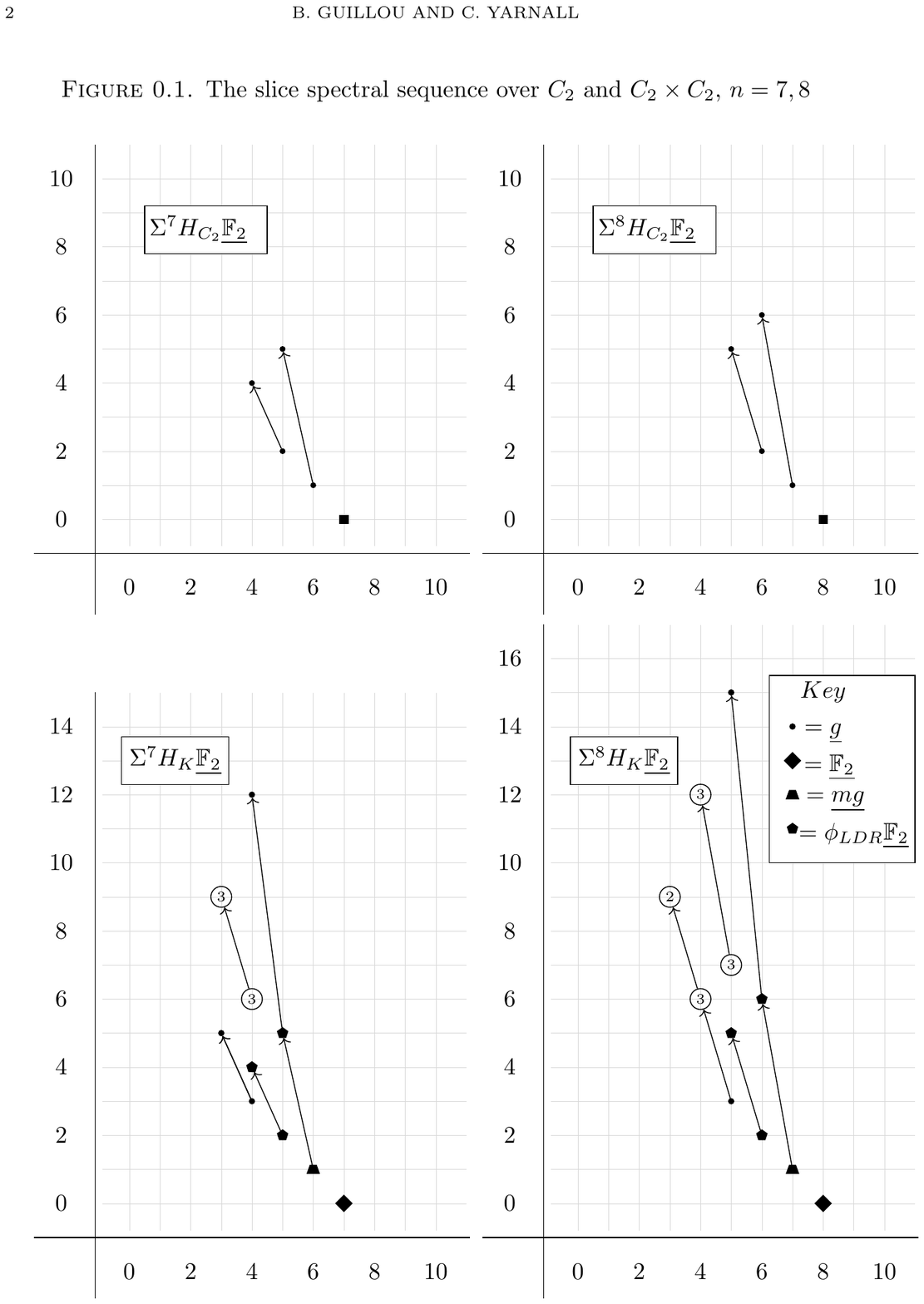}

\newpage
\begin{figure}
\caption{\label{910Fig} The slice spectral sequence over $C_2$ and $C_2\times C_2$, $n=9,10$}
\end{figure}
\includegraphics{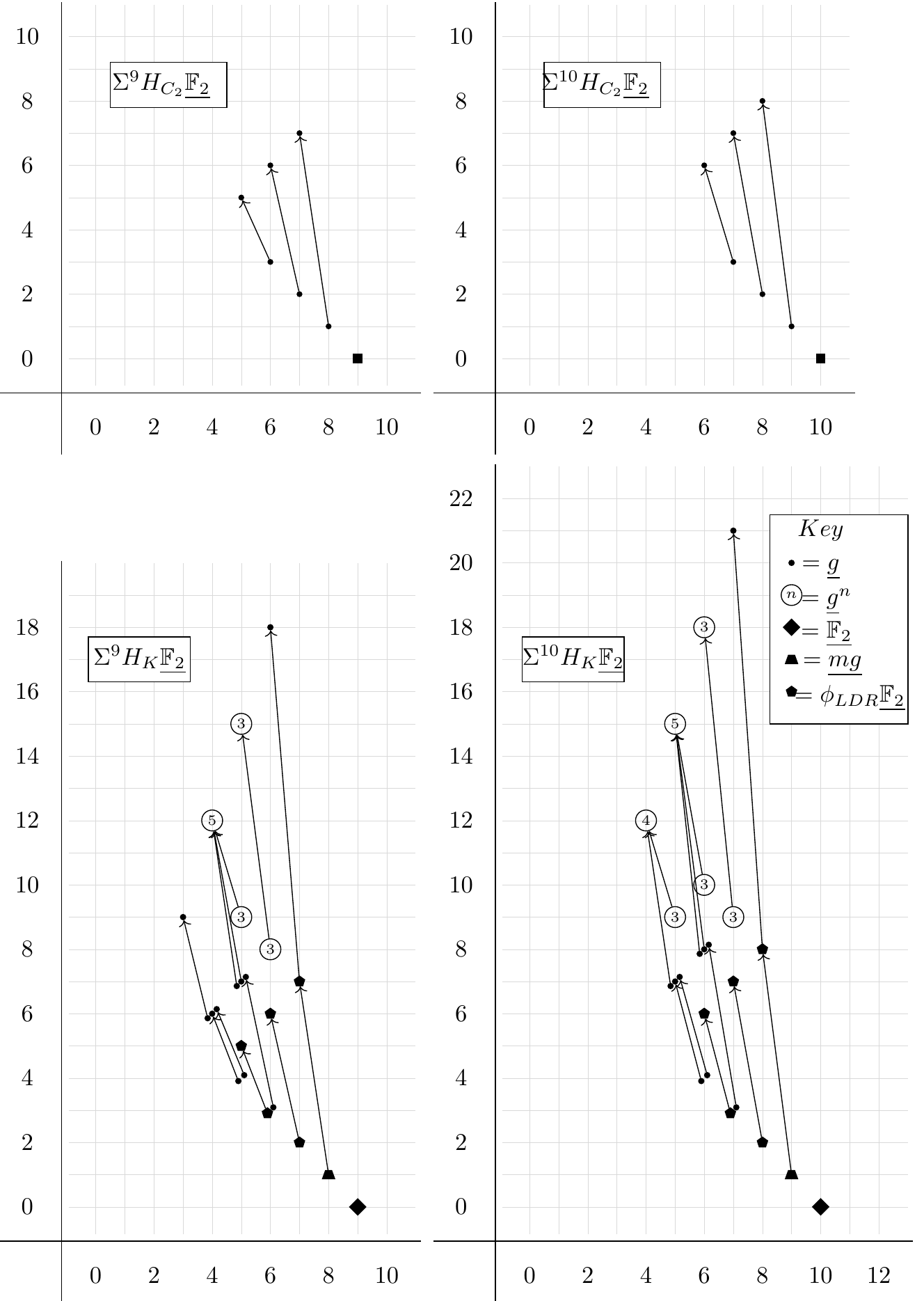}

\newpage
\begin{figure}
\caption{\label{1112Fig} The slice spectral sequence over $C_2\times C_2$, $n=11,12$}
\end{figure}
\includegraphics{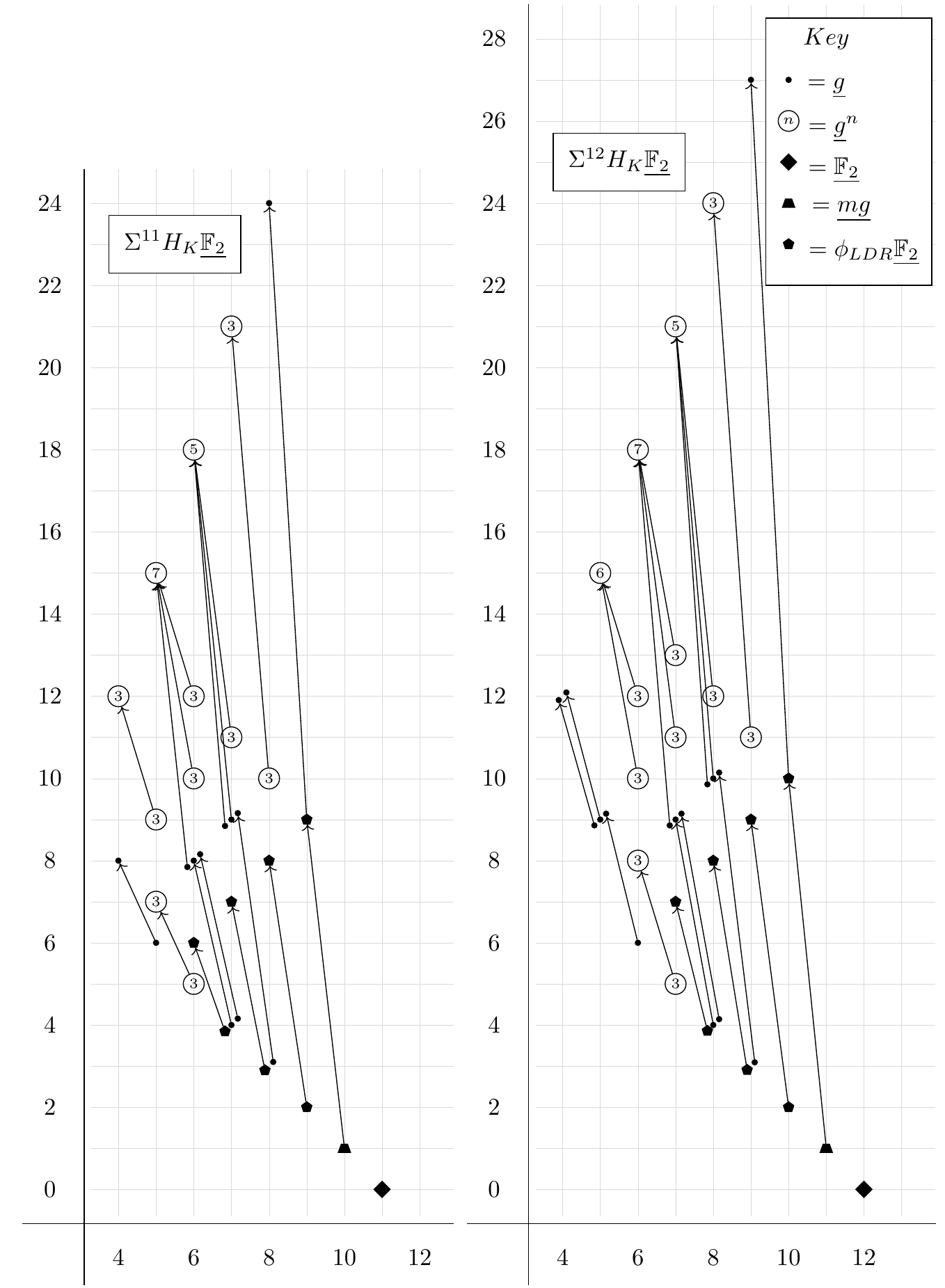}

\newpage
\begin{figure}
\caption{\label{20Fig} The slice spectral sequence over $C_2\times C_2$, $n=20$}
\end{figure}
\includegraphics{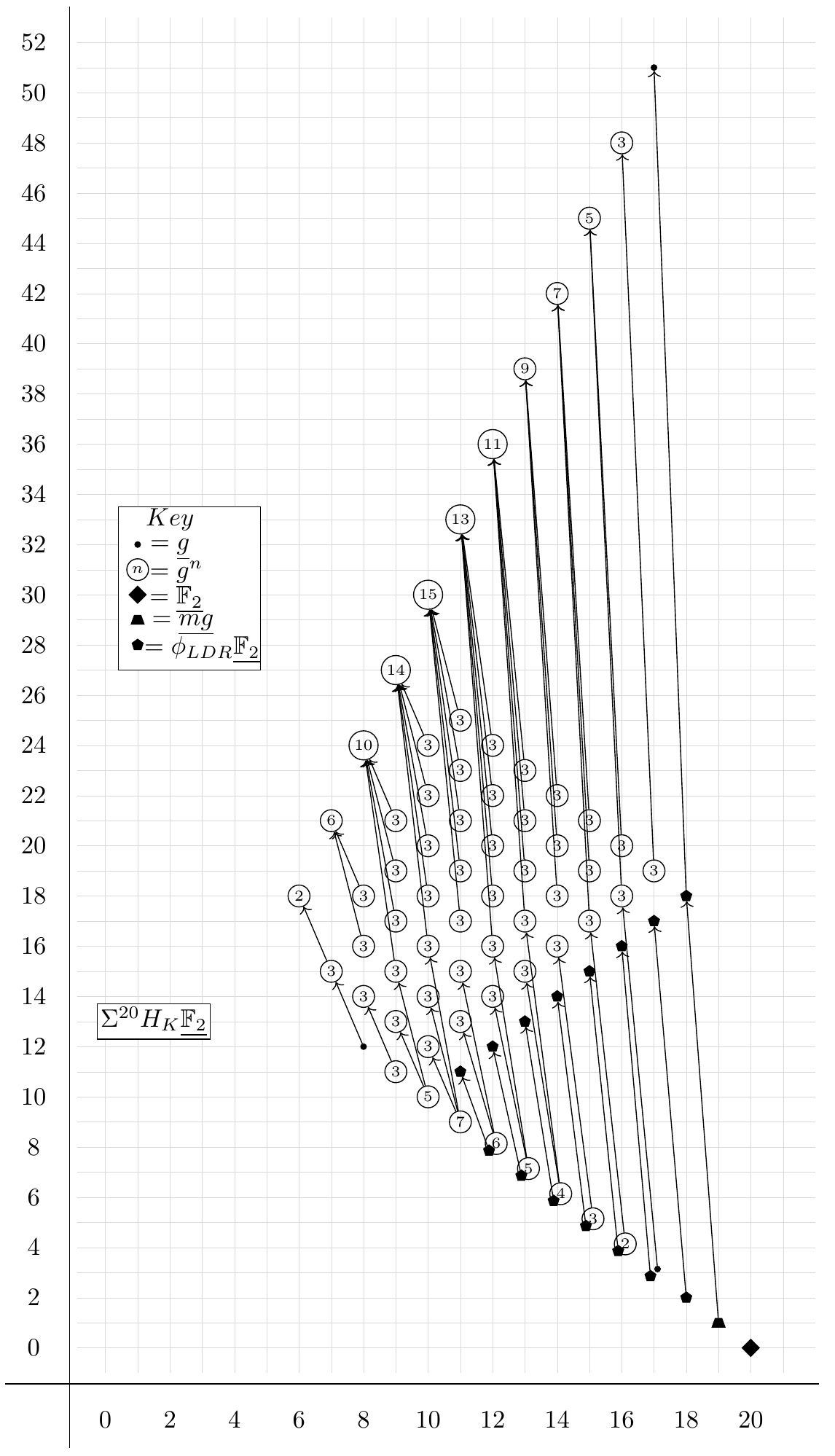}

\section*{Appendix: 
Mackey functors}\label{sec:app}

\begin{center}
\begin{tabular}{LLLl} 
\hline
\text{name} & \text{Mackey functor} & \text{description} & reference \\ \hline
\ulF & \raisebox{3.8em}{\xymatrix{
 & \F \ar[dl]_{1} \ar[d]_{1} \ar[dr]^1 & \\
 \F \ar[dr]_1& \F \ar[d]_1  & \F \ar[dl]^1,   \\ 
  & \F  &
}}   & \\  \cline{2-2}
\ulF^* & \raisebox{3.8em}{\xymatrix{
 & \F \ & \\
 \F \ar[ur]^1& \F \ar[u]^1  & \F \ar[ul]_1,   \\ 
  & \F  \ar[ul]^1  \ar[u]^1 \ar[ur]_1 &
}}  & H\ulF^* \simeq \Sigma^{4-\rho} H\ulF & \autoref{RhoDesusp} \\ \cline{2-2}
\ulg & \raisebox{3.8em}{\xymatrix{
 & \F  & \\
 0 & 0  & 0,   \\ 
  & 0  &
}} & \ulg = \phi_\Kl^*(\F)  \\  \cline{2-2}
\ulf & \raisebox{3.8em}{\xymatrix{
 & 0 & \\
 0 & 0 & 0 \\
  & \F   &
}}
&  \\  \cline{2-2}
\phi_{LDR}^* \ulf & \raisebox{3.8em}{\xymatrix{
 & 0 & \\
 \F & \F & \F \\
  & 0   &
}}
& & \autoref{InflNotn} \\  \cline{2-2}
\phi_{LDR}^* \ulF & \raisebox{3.8em}{\xymatrix{
 & \F\oplus \F\oplus \F \ar[dl]_{p_1} \ar[d]_{p_2} \ar[dr]^{p_3} & \\
 \F & \F & \F \\
  & 0   &
}}
& &  \autoref{InflNotn} \\  
\hline
\end{tabular}


\begin{tabular}{LLLl} 
\hline
\text{name} & \text{Mackey functor} & \text{description} & reference \\ \hline
 \mf{m} &
\raisebox{3.8em}{\xymatrix{
 & \F \ar[dr]^1 \ar[dl]_1 \ar[d]_1 & \\
 \F   & \F  & \F  .  \\
  & 0 &
}} & \coker \big( \, \ulf \rtarr \ulF \, \big )& \\ \cline{2-2}
\mf{mg} &
\raisebox{3.8em}{\xymatrix{
 & \F \oplus \F \ar[dr]^{p_2} \ar[dl]_{p_1} \ar[d]_{\nabla} & \\
 \F   & \F  & \F  .  \\
  & 0 &
}} & H\mf{mg} \simeq \SI^{\rho-2} H\mf{m}^* & \autoref{mProp} \\ \cline{2-2}
\mf{w} &
\raisebox{3.8em}{\xymatrix{
 & 0  & \\
 \F  \ar[dr]_1 & \F \ar[d]^1  & \F \ar[dl]^1  .  \\
  & \F &
}} & \ker \big(\, \ulF \rtarr \ulg\, \big) &
\\  \cline{2-2}
\mf{W} & 
\raisebox{3.8em}{\xymatrix{
 & \F \oplus \F \oplus \F & \\
 \F  \ar[ur]^{\iota_1} \ar[dr]_1 & \F \ar[u]^{\iota_2} \ar[d]^1  & \F \ar[ul] _{\iota_3} \ar[dl]^1.  \\
  &  \F   &
}} & \\
\hline
\end{tabular}
\end{center}

\end{document}